\pdfoutput=1
\documentclass[mnsc]{informs_modified}

\OneAndAHalfSpacedXI

\usepackage{natbib}
 \bibpunct[, ]{(}{)}{,}{a}{}{,}%
 \def\bibfont{\small}%
\usepackage{color}
\newcommand{\x}{\textcolor{black}}
\newcommand{\va}{\textcolor{black}}
\newcommand{\bE}{\mathbf{E}}
\newcommand{\mujL}{{\mu_j^L}}
\newcommand{\yjS}{{Y_j^S}}
\newcommand{\yjL}{{Y_j^L}}
\newcommand{\wkS}{{W^S_k}}
\newcommand{\ratioL}{{\text{ratio}_j^L}}
\newcommand{\ratioS}{{\text{ratio}_j^S}}
\newcommand{\ratioall}{{\text{ratio}_j^\text{all}}}
\newcommand{\OPTI}{\text{OPT}(I)}
\newcommand{\ALGI}{\text{ALG}(I)}

\usepackage{array}
\usepackage{multirow}
\newcolumntype{L}{>{\centering\arraybackslash}m{3cm}}
\TheoremsNumberedThrough     
\ECRepeatTheorems

\EquationsNumberedThrough    

\MANUSCRIPTNO{MS-0001-1922.65}

\begin{document}


\RUNAUTHOR{Stein, Truong, and Wang}

\RUNTITLE{Advance Service Reservations with Heterogeneous Customers}

\TITLE{Advance Service Reservations with Heterogeneous Customers}

\ARTICLEAUTHORS{%
\AUTHOR{Clifford Stein, Van-Anh Truong, Xinshang Wang} \AFF{Department
of Industrial Engineering and Operations Research, Columbia
University, New York, NY 10027, USA, \EMAIL{cliff@ieor.columbia.edu, vatruong@ieor.columbia.edu, xw2230@columbia.edu}
\URL{}}
}

\ABSTRACT{We study a fundamental model of resource allocation in which
  a finite number of resources must be assigned in an online manner to
  a heterogeneous  stream of customers. The  customers arrive randomly
  over time  according to  known stochastic processes.   Each customer
  requires a specific amount of capacity and has a specific preference
  for each  of the resources,  with some resources being  feasible for
  the  customer  and  some  not.   The system  must  find  a  feasible
  assignment  of  each customer  to  a  resource  or must  reject  the
  customer.   The  aim is  to  maximize  the  total expected  capacity
  utilization  of the  resources  over the  horizon.   This model  has
  application   in  services,   freight  transportation,   and  online
  advertising.  We present  online algorithms with bounded competitive
  ratios  relative to  an  optimal offline  algorithm  that knows  all
  stochastic  information.   Our  algorithms  perform  extremely  well
  compared to common heuristics, as demonstrated on a real data set
  from a large hospital system in New York City.}

\KEYWORDS{Analysis of algorithms, Approximations/heuristic, Cost analysis} 

\maketitle

%


\section{Introduction}
We study a fundamental model of resource allocation in which a finite
number of resources must be assigned in an online manner to a
heterogeneous stream of customers. The customers arrive randomly over
time according to known stochastic processes.  Each customer requires
a specific amount of capacity and has a specific preference for each
of the resources, with some resources being feasible for the customer
and some not.  The system must find a feasible assignment of each
customer to a resource or must reject the customer.  The aim is to
maximize the total expected capacity utilization of the resources over
the time horizon.

This model has application in multiple areas, including services,
online advertising, and freight transportation.  We now explain a few
of the applications.
\begin{description} 
\item[\textbf{Service Reservation.}] In services such as healthcare,
  the resources can correspond to service sessions.  For example, a
  resource might be a Monday afternoon session from 1 to 5 PM with
  Dr. Smith. The customers are patients who arrive to book
  appointments over time.  Based on a patient's urgency, type of
  visit, arrival time, and preferences, the patient might require a
  specific length of visit and might be preferably assigned only to a
  subset of sessions.  Upon the arrival of a patient, the
  system has to reserve a part of a session for the patient.  This
  appointment decision typically takes place immediately.
If an appointment cannot be found, the system must reject
  the patient.

\item[\textbf{Generalized Adwords.}] In online advertizing, the resources correspond to advertisers.  The capacity of each resource corresponds to the budget of the corresponding advertiser.  Ad impressions arrive randomly over time. 
Each impression, depending on its characteristics, commands a known non-negative bid from each of the advertisers.  When an impression occurs, the ad platform must allocate it to an advertiser for use.  The ad platform earns the bid, and the budget of the advertiser is depleted by the same amount.  The aim of the ad platform is to maximize the expected revenue earned. Our model is more general than adwords models, as we allow bids to have arbitrary sizes, whereas adwords model tend to assume that bid sizes must be very small relative to the budgets, or that  each bid must be truncated by the remaining budget \citep{mehta2012online}.

\item[\textbf{Freight Allocation.}]
Freight carriers such as motor carriers, railroad companies, and shipping companies have fleets of containers that can be deployed to move loads from specific origins to destinations.  The assignment of containers to routes are tactical decisions that are performed on a larger time scale.  Suppose that we focus on a single route.  Each container, with its specific departure and arrival time, corresponds to a resource.  Customer demands for the route arrive randomly over time.  Each demand unit has a specific size and delivery time line.  As each demand unit arrives, the operational decision is how to assign the demand unit to a specific container in the fleet \citep{spivey2004dynamic}. This assignment generates a quoted time of delivery for the customer, reduces the available capacity in the container, and earns the system an amount that can be roughly proportional to the amount of capacity consumed.   
\end{description}

Our model captures most, if not all, of the features of the above applications. Specifically, we consider a continuous-time planning horizon.  There are $m$ resources with known capacities.  There are $n$ customer types. Each customer type is associated with a known stochastic arrival process.   Each customer can be assigned to a known subset of the resources, and consumes a known amount of each resource that it is assigned to.  The system aims to assign customers to resources immediately and irrevocably as they arrive in order to maximize the total expected amount of resources used.

A salient feature of our model is that the resources may be
perishable.  This feature makes the model especially appropriate for
service applications.  More specifically, each resource may be
associated with a known expiry date that falls within or beyond the
horizon.  The way we capture an expiry date for a specific resource is
to make that resource infeasible for all customer types that arrive
after its expiry date.  That is, to capture the perishability of
resources, we equivalently force the composition of customer types
that arrive over time to change over time. For this reason,
the non-stationarity of arrivals in our model is of especial
importance.

Another significant advantage of non-stationary arrivals is the ability to better capture real applications.  In real applications, demands can be highly non-stationary, changing with the time of day, time of week, seasons and longer-term trends \citep{huh2012multi}.  \citet{kim2014call} have shown, for example, that call-center and hospital demands are well-modeled by non-homogeneous Poisson processes.  For a problem that essentially aims to match demand with supply over time, capturing this non-stationarity in demand arrivals can lead to significant improvements in performance over stationary models.

Our basic model can be adapted as needed to fit various applications.
In this paper, we will focus on solution methods for the core model.
It is easy to see that the associated dynamic stochastic optimization
problem is intractable to solve with dynamic prog	ramming.  The state
of the system grows exponentially with the length of the horizon.
Therefore, we aim to develop near-optimal policies that are robust and
easy to compute.  We will study an online version of the problem. A
problem is {\it online} if at all points in time, the algorithm has to
make adaptive decisions based only on past and current information.  In
contrast, an {\it offline} algorithm knows all future (stochastic)
information up-front.  We will use {\it competitive analysis} to
evaluate our algorithms \citep{borodin2005online}.  We will consider
the relative expected performance between an online algorithm and an
optimal offline algorithm.  We define the minimum ratio between the
benefit achieved under the online algorithm and that under the optimal
offline algorithm as {\it competitive ratio} for that online
algorithm. An algorithm with a competitive ratio of $\alpha$ is {\it
  $\alpha$-competitive}.

We propose $0.321$-competitive online algorithms.  Further, we show that an upper bound on the competitive ratio of any algorithm is 1/2.   Ours are the first algorithms with performance guarantees for the advance reservation of service with heterogeneous customer needs and preferences.  They are also the first algorithms with constant competitive ratios for the adwords problem without any assumption on the bid size and on the stationarity of the arrival process.  Despite the conservative performance characterization, we show that our algorithms perform extremely well compared to common heuristics as demonstrated on a real data set from a large hospital system in New York City.

\section{Literature Review}
Our model is related to many streams of literature, the closest of which are the adwords, the dynamic knapsack, and the appointment-scheduling literature.

\subsection{Adwords problems}
Our model generalizes adwords problems.  Considerable work has been done in this area.  If  each bid is truncated by the remaining budget, it was shown by \cite{mehta2012online} that a greedy algorithm achieves a worst-case competitive ratio of $1/2$ in the adversarial-demand model.  For adwords models in which demands arrive in random orders and bids are small, \cite{goel2008online} prove that a greedy algorithm achieves a worst-case ratio of $1-1/e$.  \citet{mirrokni2012simultaneous} later improve this ratio to 0.76.  If demands are i.i.d., but bids are not necessarily small, \citet{devanur2011near} show that a greedy algorithm achieves the worst-case ratio of $1-1/e$.  Later, \citet{devanur2012asymptotically} show that under stochastic demands, if the bid to budget ratio is at most $1/d$, $d\geq 2$, and if bids can be truncated, then there is an algorithm that achieves a worst-case ratio of $1 - 1/\sqrt{2\pi d}$.  If the bid to budget ratios at most $\epsilon^2$, then the algorithm achieves a worst-case ratio of $1-O(\epsilon)$.  Finally, no algorithm can achieve a worst-case ratio that is better than $1-o(1/\sqrt{d})$ when the bid to budget ratios are as large as $1/d$.  The main difference between our work and this literature is that we do not make the assumption of truncated bids, small bids, or i.i.d. demand.  Furthermore, we study the ratio of expected performance between the online and optimal offline algorithm, rather than the worst-case ratio.

\begin{table}
\begin{center}
\caption{Results on adwords models.}
\label{tab:adwords}
\begin{tabular}{|c|c|c|L|}
\hline
Reference & Lower bound achieved & Assumption\\
\hline
\hline
\cite{mehta2012online} & 0.5 & adversarial demand, truncated bids\\
\hline
\multirow{ 2}{*}{Our work} & 0.321 & stochastic demand\\\cline{2-3}
 & $1 - 1/\sqrt{2\pi d} + O(1/d)$ & stochastic demand, bid to budget ratio at most $1/d$\\ 
\hline
\cite{goel2008online} & $1-1/e\approx 0.63$ & randomly ordered demand, small bids\\
\hline
\cite{mirrokni2012simultaneous} & 0.76 & randomly ordered demand, small bids\\
\hline
\cite{devanur2012asymptotically}&$ 1 - 1/\sqrt{2\pi d}$ & \begin{tabular}{c}stochastic demands,\\ bid to budget ratio at most $1/d$, $d\geq 2$,\\ truncated bids\end{tabular}\\
\hline
\cite{devanur2011near} & $1-1/e\approx 0.63$ & i.i.d. demand, truncated bids\\
\hline
\end{tabular}
\end{center}
\end{table}

\subsection{Dynamic knapsack problems}
Our problem is related to multi-constrained dynamic knapsack problems (MKP).  In these problems, a set of randomly arriving items must be packed into one or more knapsacks, respecting the capacity constraints of the knapsacks.  The goal is to maximize the value of the items packed.  Note that our problem is different from these dynamic multi-knapsack problems.  In our problem, each customer can be satisfied using one of a subset of resources, rather than any resource, due to preferences, urgency, priorities, etc.  These feasibility constraints must be accounted for in the assignment decision.  In contrast, a knapsack problems, an object can be placed into any knapsack, as long as the capacity constraints are satisfied.

Dynamic-programming characterizations have been studied in the case of one knapsack \citep{papastavrou1996dynamic, kleywegt1998dynamic,van2000finite, lin2008stochastic, chen2014adaptive}.  Some results generalize to the MKP but these results are not sufficient to yield provable approximations \citep{van2000finite}.  Many authors have studied online algorithms for the MKP. It is shown in \cite{marchetti1995stochastic} that no online algorithm for MKP exists with a constant worst-case competitive ratio.  Therefore, \cite{marchetti1995stochastic} and \cite{lueker1998average} study algorithms with bounded additive differences away from the offline optimal. Finally, \citet{chakrabarty2013online} design an algorithm with a bounded worst-case competitive ratio, assuming that the size of each item is very small relative to the capacity, and the value-to-size ratio of each item is upper and lower bounded by two constants.  Our model is different from the MKP model because our resources are not interchangeable, as customer preferences for them might be different.  Our approach also differs from existing MKP approaches in that we seek to bound the ratio of expected performance between the online and optimal offline algorithm, rather than the worst-case ratio.

\subsection{Online resource-allocation problems} \label{sec:onlineProblems}
Our model falls within the literature on online resource allocation.  Adwords and dynamic knapsack problems are subclasses of this literature.  Although this is a vast literature, it can be roughly divided into several streams based on the assumptions made about the model.  The first stream is focused on designing algorithms for problems in which demands arrive in adversarial fashion \citep{karp1990optimal,kalyanasundaram1996optimal,aggarwal2011online,devanur2013randomized}.  The second stream is focused on problems in which demands arrive as a random permutation of a known sequence \citep{goel2008online,agrawal2009dynamic,devanur2009adwords,mahdian2011online,karande2011online,bhalgat2012online}.  The third stream is focused on problems in which demands are drawn i.i.d. from an unknown distribution \citep{ghosh2009bidding,devanur2011near,balseiro2014yield}.   
 The fourth stream is focused on problems in which demands are drawn i.i.d. from an known distribution \citep{feldman2009online,agrawal2009dynamic,feldman2010online,vee2010optimal,jaillet2012near,manshadi2012online,jaillet2013online}.  The fifth stream is focused on resource allocation under the small-bid or truncated-bid assumption \citep{mehta2007adwords,buchbinder2007online,jaillet2011online,devanur2012asymptotically}.  Very few papers focus on models with non-stationary stochastic demand as we do, and in these cases, they \x{either assume that bids can be truncated \citep{ciocan2012dynamic}}, or assume that the resource consumption is constant for all demand units, i.e., the problems are matching problems \citep{alaei2011adcell,wangTB2015,gallegoLTW2015}.

\subsection{Appointment-scheduling problems}

Our work is related to the literature on appointment scheduling. 	Most relevant is the stream of literature focusing on how to assign future appointments to patients. This paradigm is called {\it advance scheduling}.  In the literature of advance scheduling, \citet{Truong2014Advance} studies the analytical properties of a two-class advance-scheduling model and gives efficient methods  for computing an optimal scheduling policy.  \citet{gocgun2012lagrangian}, \citet{patrick2008dynamic}, \citet{feldman2014appointment}, and \citet{gupta2008revenue} study structural properties of optimal policies and propose heuristics for several related models, although they do not investigate the theoretical performance of these heuristics.  \citet{wangTB2015} propose online algorithms with constant competitive ratios for advance scheduling with multiple patient types and with patient preferences.  Their model is very close to ours in that the demand processes are allowed to be known non-stationary stochastic processes.  They also define the same notion of competitive ratio.  However, their model is considerably easier, since they assume that each customer demands a unit amount of capacity.  In contrast, we assume that customers have heterogeneous capacity requirements.  As we shall show, much of the effort in our algorithms and their analysis is directed towards taking care of these differences in capacity requirement.

The paper is organized as follows.  We specify the model and performance metric in Section \ref{sec:model}. In Section \ref{sec:bound}, we prove that $0.5$ is an upper bound on the competitive ratio of any online algorithm for this problem.  We derive an upper bound on the optimal offline objective in Section \ref{sec:upperBound}. In Section \ref{sec:warmup}, we design a basic online algorithm with a competitive ratio of $0.5 (1-1/e) \approx 0.316$, which serves to illustrate our key ideas.  In Section \ref{sec:improve}, we refine the algorithm to employ resource sharing in order to obtain an improved competitive ratio of $0.321$, as well as an improved empirical performance.  In Section \ref{sec:numerical}, we compare the empirical performance of our algorithms against two commonly used heuristics by simulating the algorithms on appointment-scheduling data obtained from a large hospital system in New York City.

\section{Model and Performance Metric}\label{sec:model}
We use $[n]$ to denote the set $\{1,2,\ldots,n\}$ and   consider a continuous
horizon $[0,T]$.  There are $m$ resources and $n$ customer types.
Resource $j \in [m]$ has
capacity $c_j \in \mathbb{R}_+$.  
Customers of type $i \in [n]$ arrive according to a non-homogeneous
Poisson process with rate $\lambda_i(t)$, for $t\in [0,T]$. Let $\Lambda_i = \int_0^T \lambda_i(t) dt$ be the expected total number of arrivals of type-$i$ customers. The
arrival rates of all the customer types are known.  When a customer
arrives, one of the $m$ resources needs to be immediately allocated to
the customer, or the customer must be rejected.  If resource $j$ is
allocated to a customer of type $i$, exactly $u_{ij}$ units of
resource $j$ must be provided. We assume that the $u_{ij}
\in [0,c_j]$, $\forall i \in [n]$ and $j \in [m]$, are known.  The
reward earned for the assignment of customer type $i$ to resource $j$
is also $u_{ij}$.  The objective is to maximize the total expected reward
over the horizon, which equivalently maximizes total resource
utilization.

Let $I$ be a sample path of customer arrivals over the entire horizon. We say that an algorithm is \emph{offline} if it knows $I$ at time $0$. An algorithm is \emph{online} if at any time $t$, it only knows future arrival rates and the realization of all the arrivals prior to $t$.

Let $\ALGI$ be the total amount of resources allocated by an online algorithm $ALG$. Let $\OPTI$ be the total amount of resources allocated by an optimal offline algorithm $OPT$. We define the competitive ratio of $ALG$ to be
\begin{equation} \text{Competitive Ratio of $ALG$} = \frac{\mathbf{E}[\ALGI]}{\mathbf{E}[\OPTI]},\label{eq:cr}\end{equation}
where the expectation is taken over all sample paths $I$ and over the random realizations of the online algorithm.  Our definition of competitive ratio follows previous works including \citet{feldman2009online}, \cite{jaillet2013online}, and \citet{wangTB2015}.  It is less conservative than the worst-case ratio $\min_I \frac{\ALGI}{\OPTI}$ that has been more commonly used for online algorithms.

\section{Upper Bound on the Competitive Ratio}\label{sec:bound}
\x{In this section, we prove upper bounds on competitive ratios that can be achieved by any online algorithm. }
\x{We first prove a uniform upper bound on the competitive ratio.} 

\begin{proposition}
The competitive ratio of any online algorithm is at most $0.5$.
\end{proposition}
\proof{Proof.} 
Consider an input with two customer types and a single resource. Assume that the horizon is $[0,1]$.  The capacity of the resource is $c_1 = 1$.
\begin{itemize}
\item Type-1 customers have a very large arrival rate in time $[0,0.5]$, but their arrival rate is $0$ after time $0.5$. In particular, $\Lambda_1 = \int_0^{0.5} \lambda_1(t) dt \gg 1$, so that we can ignore the event that no type-1 customer arrives. Their utilization for the single resource is $u_{11} = \epsilon / \Lambda_1$ for some very small value $\epsilon$.
 \item Type-2 customers arrive in time $(0.5,1]$. They have a very small arrival rate $\Lambda_2 = \int_{0.5}^1 \lambda_2(t) dt = \epsilon$. Their utilization for the resource is $u_{21} = c_1 = 1$.
\end{itemize}

Since customers of type 2 request the entire resource, the offline algorithm will allocate the resource to a type-2 customer if there is one. The probability that at least one type-2 customer arrives is $1 - e^{-\Lambda_2} = \Lambda_2 + o(\Lambda_2^2) = \epsilon + o(\epsilon^2)$. With probability $1 - o(\Lambda_2) = 1 - o(\epsilon)$, no type-2 customer will arrive, in which case the optimal offline algorithm will accept as many type-1 customers as possible. The expected total utilization of all type-1 customers is $u_{11} \cdot \Lambda_1 = \epsilon$. Suppose $\epsilon \ll c_1 = 1$. Then all type-1 customer can be accepted. In sum, the expected amount of resource allocated by an optimal offline algorithm is \begin{align*}
&   1\cdot (\epsilon + o(\epsilon^2))  + \epsilon \cdot (1 - o(\epsilon))\\
= & 2 \epsilon +  o(\epsilon^2).
\end{align*}

The decision of an online algorithm is whether to accept type-1 customers during time $[0,0.5]$. If it does accept type-1 customers, the online algorithm earns  $u_{11} \cdot \Lambda_1 = \epsilon$ in expectation. Otherwise, with probability $\Lambda_2 + o(\Lambda_2^2)$ it earns $u_{21}$, which is $u_{21} (\Lambda_2 + o(\Lambda_2^2)) = \epsilon + o(\epsilon^2)$ in expectation. In sum, an online algorithm cannot allocate more than $\epsilon + o(\epsilon^2)$ in expectation. Thus, an upper bound on the competitive ratio is
\[(\epsilon + o(\epsilon^2)) / (2\epsilon + o(\epsilon^2)),\]
which tends to 0.5 in the limit as $\epsilon \to 0$.
\halmos
\endproof

\x{Next, we prove upper bounds on competitive ratios for special cases in which the utilization $u_{ij}$ for each resource is bounded away from $c_j$. Specifically, suppose there is some integer $d \geq 2$ for which $u_{ij} \leq c_j / d$  for all $i\in[n], j \in[m]$. In Proposition \ref{prop:upperboundd}, we prove an upper bound that depends on any finite $d$. In Proposition \ref{prop:upperboundd2}, we prove an even tighter upper bound for the asymptotic regime $d \to \infty$. } We introduce the following technical lemma for proving the parameter-dependent bounds on competitive ratios.

\begin{lemma}\label{lm:Revision2UpperBound}
If $N$ is a Poisson random variable and 
\[ \sum_{i=1}^{d-1} P(N=i) \frac{i}{d} + P(N \geq d) = 1-\frac{1}{2d}\]
for some integer $d \geq 2$, then
\[ P(N \geq d) \geq \frac{1}{2}.\]
\end{lemma}
\proof{Proof.}
It is easy to deduce that
\begin{align}
& \sum_{i=1}^{d-1} P(N=i) \frac{i}{d} + P(N \geq d) \label{eq:Revision2UpperBoundLemma1}\\
\leq & \sum_{i=0}^{d-1} P(N=i) \frac{d-1}{d} + P(N \geq d) \nonumber \\
= & 1 - (1 - P(N \geq d))\frac{1}{d}. \label{eq:Revision2UpperBoundLemma2}
\end{align}
Thus, (\ref{eq:Revision2UpperBoundLemma2}) is an upper bound on (\ref{eq:Revision2UpperBoundLemma1}). In order to satisfy
 (\ref{eq:Revision2UpperBoundLemma1})$=1 - \frac{1}{2d}$, we must have (\ref{eq:Revision2UpperBoundLemma2})$\geq 1 - \frac{1}{2d}$. That is,
\[ 1 - (1 - P(N \geq d))\frac{1}{d} \geq 1 - \frac{1}{2d}\]
\[ \Longrightarrow P(N \geq d) \geq \frac{1}{2}.\]
\halmos
\endproof

\begin{proposition}\label{prop:upperboundd}
For any given integer $d \geq 2$, if $u_{ij} \leq c_j / d$ for all $i \in [n], j \in [m]$, then the competitive ratio of any online algorithm is at most $\frac{4d-2}{4d-1} \leq 1 - \frac{1}{4d}$.
\end{proposition}
\proof{Proof.}
Consider a single resource with capacity $c_1 = 1$. Consider two demand types that arrive over horizon $[0,1]$:
\begin{itemize}
\item Type-1 customers only arrive during the first half of the horizon, i.e., $[0,0.5]$. During the first half of horizon, their arrival rate is huge $\Lambda_1 = \int_0^{0.5} \lambda_1(t) dt \gg 2d$, so that we can ignore the event that fewer than $2d$ type-1 customers arrive. Assume $u_{11} = \frac{1}{2d} + \epsilon$ for some infinitesimally small $\epsilon$.
\item Type-2 customers  arrive only during the second half of the horizon, i.e., $(0.5, 1]$. Let $N$ be the total number of arrivals of type-2 customers. Assume that $u_{21} = 1/d$.
\end{itemize} 

The optimal online algorithm has only two choices: (i) accept at least one type-1 customer, in which case the total reward is at most $1-\frac{1}{2d} + O(\epsilon)$ (because the $\epsilon$ allocation forbids the last $\frac{1}{2d}$ unit of the resource from being taken); (ii) do not accept any type-1 customer, in which case the total reward is only collected from type-2 customers
\[ \sum_{i=1}^{d-1} P(N = i) i \cdot u_{21} + P(N \geq d) d \cdot u_{21} = \sum_{i=1}^{d-1} P(N=i) \frac{i}{d} + P(N\geq d) .\]

Assume that the distribution of $N$ is such that 
\[\sum_{i=1}^{d-1} P(N=i) \frac{i}{d} + P(N\geq d)  = 1-\frac{1}{2d}.\]
Then the expected total reward earned by the optimal online algorithm is at most $1-\frac{1}{2d} + O(\epsilon)$.

The optimal offline algorithm also chooses from the above two options, but now first observes the value of $N$. If $N \geq d$, the optimal offline algorithm accepts $d$ customers of type 2. Otherwise, the optimal offline algorithm accepts $2d-1$ type-1 customers to achieve total reward $1 - \frac{1}{2d} + O(\epsilon)$. In sum the expected total reward of the offline algorithm is
\[ P(N \geq d)  + P(N < d) (1 - \frac{1}{2d} + O(\epsilon)).\]

By Lemma \ref{lm:Revision2UpperBound}, we have $P(N \geq d) \geq \frac{1}{2}$. Then the expected total reward of the offline algorithm can be bounded by
\begin{align*}
&  P(N \geq d)  + P(N < d) (1 - \frac{1}{2d} + O(\epsilon))\\
= & P(N \geq d)  + P(N < d)(1 - \frac{1}{2d}) + O(\epsilon)\\
= & P(N \geq d)  + ( 1- P(N \geq d))(1 - \frac{1}{2d}) + O(\epsilon)\\
= & P(N \geq d)\frac{1}{2d}  +  1 - \frac{1}{2d} + O(\epsilon)\\
\geq & \frac{1}{2} \cdot \frac{1}{2d}  +  1 - \frac{1}{2d} + O(\epsilon)\\
= & 1 - \frac{1}{4d} + O(\epsilon).
\end{align*}

In sum, the competitive ratio of the optimal algorithm is at most
\[ \frac{1 - \frac{1}{2d} + O(\epsilon)}{1 - \frac{1}{4d} + O(\epsilon)} = \frac{4d -2 + O(\epsilon)}{4d -1 + O(\epsilon)}.\]

When $\epsilon$ tends to $0$, the ratio becomes 
\[ \frac{4d-2}{4d-1} \leq 1 - \frac{1}{4d}.\]
\halmos
\endproof

\begin{proposition}\label{prop:upperboundd2}
For large positive integer $d$, if $u_{ij} \leq c_j / d$ for all $i \in [n], j \in [m]$, then the competitive ratio of any online algorithm is at most $1 - \frac{1}{2\sqrt{\pi d}}+ o(1/\sqrt{d})$.
\end{proposition}
\proof{Proof.}
We construct a special case of our model as follows. Let $c_j = 1$ for all $j\in[m]$. There are $n = m+1$ customer types.
\begin{itemize}
\item For each customer type $i=1,2,...,m$, $u_{ii}=1/d$ and $u_{ij}=0$ for all $j\not= i$.
\item For customer type $m+1$, $u_{m+1,j}=\epsilon$ for all $j \in [m]$. We will let $\epsilon$ tend to $0$.
\item For each customer type $i=1,2,...,m$,
\[ \lambda_i(t) = \left\{ \begin{array}{ll} 0, & \forall t \in [0,T/2)\\ d/T , & \forall t \in [T/2,T] . \end{array} \right.\]
As a result, $\Lambda_i = \int_0^T \lambda_i(t) dt = d/2$ for all $i=1,2,...,m$.
\item For customer type $m+1$, 
\[ \lambda_{m+1}(t) = \left\{ \begin{array}{ll} \frac{m}{\epsilon T}, & \forall t \in [0,T/2]\\ 0 , & \forall t \in (T/2,T] . \end{array} \right.\]
As a result, $\Lambda_{m+1} =\int_0^T \lambda_{m+1}(t) dt =  \frac{m}{2\epsilon}$.
\end{itemize}

We first analyze the optimal online algorithm. For each resource $i \in [m]$, in the second half of the horizon, customers of type $i$ will request $\Lambda_i u_{ii} = 0.5$ amount of resource $i$, i.e., half of the resource, in expectation. On the other hand, in the first half of the horizon, customers of type $m+1$ can totally take $\Lambda_{m+1} u_{m+1,j} = 0.5 m$ units of any resource in expectation. 
We let $\epsilon$ tend to $0$, so customers of type $m+1$ will request almost exactly $0.5m$ units of any resource. 

It is easy to see that the optimal online strategy is to give type-$(m+1)$ customers $0.5^-$ unit of each of the $m$ resources, where $0.5^-$ means infinitesimally approaching $0.5$ from below as $\epsilon$ tends to $0$. This is because of symmetry, or more rigorously, because the marginal reward of adding $1/d$ unit of capacity to each resource is decreasing in the amount of the resource that is remaining at time $T/2$.

Starting from time $T/2$, the optimal online algorithm assigns the remaining $0.5^+$ unit of each resource $i$ to customers of type $i$. Let $N_i$ be the number of arrivals of type-$i$ customers, for $i=1,2,...,m$. We must have $\bE[N_i] = \Lambda_i = d/2$. The  total expected amount of resource $i$ assigned to customers by the end is
\[ 0.5^- +  \bE[ \min(0.5 , N_i u_{ii})],\]
where $0.5^-$ is the amount assigned to customers of type $m+1$, and $\bE[ \min(0.5 , N_i u_{ii})]$ is the expected amount assigned to customers of type $i$.

When $d$ is large, we apply the central limit theorem so that $N_i u_{ii}$ is approximated by a normal distribution with mean $0.5$ and variance $\frac{1}{2d}$. Then
\[\bE[ \min(0.5 , N_i u_{ii})] = 0.5 -  \frac{1}{\sqrt{2d}} \int_0^\infty \phi(x) x dx  + o(1/\sqrt{d})=0.5 - \frac{1}{2\sqrt{\pi d}} + o(1/\sqrt{d}), \]
where $\phi(\cdot)$ is the standard normal pdf.

In sum, when $d$ is large so that the distribution of $N_i$ can be approximated by the central limit theorem, the optimal online algorithm earns expected reward
\[ 0.5^- + 0.5 - \frac{1}{2\sqrt{\pi d}} + o(1/\sqrt{d})= 1 - \frac{1}{2\sqrt{\pi d}}+ o(1/\sqrt{d})\]
from each resource $i$, and earns total expected reward
\[ m(1 - \frac{1}{2\sqrt{\pi d}})+ o(m/\sqrt{d})\]
from all the $m$ resources.

Now we analyze the optimal offline algorithm. The optimal offline algorithm will first fill  each resource $i$ with customers of type $i$, for $i=1,2,...,m$. After that, the total remaining capacity of all the $m$ resources will be
\[ \sum_{i=1}^m \max(0,1 - N_i u_{ii}).\]

We again apply the central limit theorem so that each $N_i u_{ii}$ is approximated by a normal distribution with mean $0.5$ and variance $\frac{1}{2d}$. Then the total remaining capacity can be written as
\[ m - \sum_{i=1}^m N_i u_{ii} - o(m/\sqrt{d}),\]
where $m$ is the total initial capacity of all the $m$ resources; $N_i u_{ii}$ is the amount of resource requested by customers of type $i$; $-o(m/\sqrt{d})$ represents the loss from the approximation by the central limit theorem, plus the loss from the (ignorable) events that $N_i u_{ii} > 1$. Now $\sum_{i=1}^m N_i u_{ii}$ is approximated by a normal random variable with mean $m/2$ and variance $\frac{m}{2d}$.

Next, the optimal offline algorithm fills the remaining capacity with customers of type $m+1$, who will totally take $(m/2)^-$ units of any resource. After that, the expected remaining capacity is
\begin{align*}
& \bE[\max(0, m - \sum_{i=1}^m N_i u_{ii} - o(m/\sqrt{d}) - (m/2)^- )]\\
= & \bE[\max(0, m/2 - \sum_{i=1}^m N_i u_{ii} - o(m/\sqrt{d}))]\\
= & \sqrt{\frac{m}{2d}} \int_0^\infty \phi(x) x dx + o(m/\sqrt{d})\\
= & \frac{\sqrt{m}}{2\sqrt{\pi d}} + o(m/\sqrt{d}).
\end{align*}

The total expected reward earned by the optimal offline algorithm is 
\[m - \frac{\sqrt{m}}{2\sqrt{\pi d}} + o(m/\sqrt{d}),\]
which is the total capacity $m$ less the expected remaining capacity.

Finally, the competitive ratio of the optimal online algorithm is
\[ \frac{ m(1 - \frac{1}{2\sqrt{\pi d}})+ o(m/\sqrt{d})}{m - \frac{\sqrt{m}}{2\sqrt{\pi d}} + o(m/\sqrt{d})} = \frac{ 1 - \frac{1}{2\sqrt{\pi d}}+ o(1/\sqrt{d})}{1 - \frac{1}{2\sqrt{m \pi d}} + o(1/\sqrt{d})},\]
which tends to $1 - \frac{1}{2\sqrt{\pi d}}+ o(1/\sqrt{d})$ when $m = d$.
\halmos
\endproof



\section{Upper Bound on the Optimal Offline Objective}\label{sec:upperBound}
We derive an upper bound on the optimal offline objective, namely  $\mathbf{E}[\OPTI]$.  Since $\mathbf{E}[\OPTI]$ is very hard to analyze due to its complex offline properties, we are interested in developing an upper bound on $\mathbf{E}[\OPTI]$, which is more tractable. We will later compare the performance of our online algorithms against this upper bound, rather than directly with $\mathbf{E}[\OPTI]$.

Our upper bound can be formulated as a static LP, which allocates the expected demands $\Lambda_i$, $i \in [n]$, to the capacities $c_j$, $j \in [m]$. The decision variable $x_{ij}$ of the LP stands for the average number of customers of type $i$ to be allocated to resource $j$. The LP produces a fractional assignment.
\begin{align}\label{eq:lp}
\begin{split}
V^{LP} = \max_{x_{ij}} & \sum_{i \in [n]} \sum_{j \in [m]} x_{ij}u_{ij}\\
\text{s.t. } & \sum_{i\in [n]} x_{ij} u_{ij} \leq c_j,  \,\,\,\,\, \forall j \in [m]\\
& \sum_{j \in [m]}  x_{ij} \leq \Lambda_i, \,\,\,\,\, \forall i \in [n]\\
& x_{ij} \geq 0, \,\,\,\, \forall i \in [n], j \in [m].
\end{split}
\end{align}
By the linearity of assignment problems, it can be shown easily that
\begin{proposition}\label{prop:offline}
$V^{LP}$ is an upper bound on $\mathbf{E}[\OPTI]$.
\end{proposition}
\noindent \textbf{Proof of Proposition \ref{prop:offline}.}
\proof{Proof.}
Let $a_i(I)$ be the actual number of arrivals of type-$i$ customers in sample path $I$. Let $\tilde x(I)$ be a corresponding optimal offline (fractional) assignment. Then $\tilde x(I)$ must satisfy
\[\sum_{i\in [n]} \tilde x_{ij}(I) u_{ij} \leq c_j,  \,\,\,\,\, \forall j \in [m],\]
\[\sum_{j \in [m]}  \tilde x_{ij}(I) \leq a_i(I), \,\,\,\,\, \forall i \in [n].\]
Taking expectation on both sides, we obtain
\[\sum_{i \in [n]} \mathbf{E}[\tilde x_{ij}(I)] u_{ij} \leq c_j,  \,\,\,\,\, \forall j \in [m],\]
\[\sum_{j\in [m]}  \mathbf{E}[\tilde x_{ij}(I)] \leq \mathbf{E}[a_i(I)] = \Lambda_i, \,\,\,\,\, \forall i \in [n].\]
These inequalities imply that $\mathbf{E}[\tilde x(I)]$ is a feasible solution to (\ref{eq:lp}). Thus $V^{LP}$ must be an upper bound on $\sum_{i\in [n],j\in [m]} \mathbf{E}[\tilde x_{ij}(I)]u_{ij}$, which proves the proposition.
\halmos
\endproof

\section{Basic Online Algorithm} \label{sec:warmup}
As a warm up, we design an online algorithm which we prove to have a competitive ratio of at least $0.5 (1-1/e) \approx 0.316$. This algorithm serves to illustrate the following two main ideas, which we will later refine to obtain an improved bound.

\begin{itemize}
\item \emph{LP-based random routing.} We make use of an optimal solution $x^*$ to the static LP (\ref{eq:lp}) to route customers to resources.  Note that this solution assigns demand to supply at an aggregate level, in the expected sense. Given a solution $x^*$, for each arriving customer of type $i \in [n]$, we randomly route the customer to each candidate resource $j\in[m]$ independently with probability $x_{ij}^*/\Lambda_i$.   We say a customer is \emph{routed} to resource $j$ if resource $j$ is chosen as a candidate resource for the customer. By random routing, we can conclude that the arrival process of type-$i$ customers who are routed to resource $j$ is a non-homogeneous Poisson process with rate $\lambda_i(t) \frac{x_{ij}^*}{\Lambda_i}$, for $t \in [0,T]$.

\item \emph{Reservation by customer type.} After the random routing
  stage, we make binary admission decisions about whether to commit
  each resource $j$ to each customer $i$ who is routed to $j$. If the
  decision is `no', we reject the customer.  We make this admission
  decision as follows. For each resource $j$, we divide the candidate
  customer types who will potentially be routed to $j$ into two sets
  based on utilization $u_{ij}$. Set $L_j \subseteq [n]$ consists of
  customer types of which the utilizations $u_{ij}$ are larger than
  $c_j/2$. Mathematically,
\[ L_j = \{ i \in [n] : u_{ij} > c_j/2\}.\]
The other set $S_j = [n] - L_j$ consists of customer types with utilization $u_{ij}$ that are at most $c_j/2$.  

For each resource $j$, our algorithm chooses one set, either $S_j$ or $L_j$, whichever has the higher expected total utilization for resource $j$.  The algorithm exclusively reserves resource $j$ for customers whose types are in the chosen set.  The algorithm rejects all customer types in the complementary set.  This step is meant to resolve conflict in resource usage among different customer types by restricting use of the resource to the most promising subset of customer types.  
\end{itemize}

\noindent {\bf Large-or-Small ($LS$) Algorithm:}
\begin{enumerate}
\item (Pre-processing step) Solve the LP (\ref{eq:lp}). Let $x^*$ be an optimal solution. 
For each resource $j$, define
\[ U_j^L \equiv  \sum_{i \in L_j} x_{ij}^*u_{ij}\]
as the amount of resource $j$ allocated to customer types in $L_j$ by the static LP. Similarly, define
\[ U_j^S \equiv \sum_{i \in S_j} x_{ij}^* u_{ij}\]
as the amount of resource $j$ allocated to customer types in $S_j$ by the LP.  

\item (Reservation step) Reserve the resource $j$ for customer types in the set $L_j$ if $U_j^L \geq U_j^S$.  Otherwise, reserve resource $j$ for customer types in the set $S_j$.

\item (Random routing step) Upon an arrival of a type-$i$ customer, randomly pick a resource $j$ with probability $x_{ij}^*/\Lambda_i$. 

\item (Admission step) If the remaining capacity of resource $j$ is at least $u_{ij}$ and $i$ belongs to the set that is reserved for $j$ then accept the customer.  Otherwise, reject the customer.
\end{enumerate}

As a consequence of the random routing process, we can separate the analysis for every resource $j \in [m]$. Define
\[ U_j \equiv U_j^L + U_j^S\]
as the total amount of resource $j$ allocated by the LP. We will show that in expectation, at least
\[ \frac{1}{2}\left(1-\frac{1}{e}\right) U_j\]
units of resource $j$ will be occupied in $LS$. 

We will use the following technical lemma, which bounds the tail expectation of demands following a compound Poisson distribution.
\begin{lemma}\label{lm:bound}
Let $X_1,X_2,X_3,...$ be a sequence of i.i.d. random variables that take values from $[0,\beta]$, for some given $\beta \in [0,\frac{1}{l}]$ with $l\geq 2$ being an integer. Let $N$ be a Poisson random variable. For any given $\alpha \in [0,1]$, if 
\[\mathbf{E}\left[\sum_{k=1}^N X_k\right] = \alpha,\]
we must have
\[ \mathbf{E}\left[ \min\left( \sum_{k=1}^N X_k, 1-\beta\right)\right] \geq \frac{1-\beta}{l-1}\mathbf{E}[\min(N',l-1)],\]
where $N'$ is a Poisson random variable with mean $\alpha(l-1)/(1-\beta)$.
In particular, when $l=2$,
\[ \mathbf{E}\left[ \min\left( \sum_{k=1}^N X_k, 1-\beta\right)\right] \geq (1-\beta)\left(1 - e^{-\alpha / (1 - \beta)}\right).\]
\end{lemma}

\noindent \textbf{Proof of Lemma \ref{lm:bound}.}
\proof{Proof.}
Let $Z_1,Z_2,Z_3,...,$ be a sequence of i.i.d. random variables each following a uniform distribution over $[0,(1-\beta)/(l-1)]$. For every $k=1,2,...,$ define a function
\[ \tilde X_k(x) \equiv \frac{1-\beta}{l-1} \mathbf{1}(Z_k < x),\]
where $\mathbf{1}(\cdot)$ denotes an indicator function. 

Since $\beta \in [0,1/l]$, we must have $\beta \leq (1 - \beta)/(l-1)$. 
It is then easy to check that for any $x \in [0,\beta]$, we have 
\[ \mathbf{E}[\tilde X_k(x)] = \frac{1 - \beta}{l-1} \cdot \frac{x}{(1 - \beta)/(l-1)} = x.\]

Thus, we have for every $k=1,2,...$,
\[ \mathbf{E}[\tilde X_k(X_k) | X_k] = X_k.\]

According to Jensen's inequality, we must have
\[ \min \left(\sum_{k=1}^N X_k,1 - \beta\right)  \geq \mathbf{E}\left[\min\left(\sum_{k=1}^N \tilde X_k(X_k), 1 - \beta\right) | X_1,X_2,\ldots\right]\]
\[ \Longrightarrow \mathbf{E}\left[\min\left(\sum_{k=1}^N X_k,1 - \beta\right)\right]  \geq \mathbf{E}\left[\min\left(\sum_{k=1}^N \tilde X_k(X_k), 1 - \beta\right)\right].\]

Since $\tilde X_k(X_k)$ is either $0$ or $(1-\beta)/(l-1)$, the term $\sum_{k=1}^N \tilde X_k(X_k)$ has the same distribution as $N'(1-\beta)/(l-1)$ where $N'$ is a Poisson random variable with mean
\[ \mathbf{E}[N'] = \frac{l-1}{1 - \beta}\mathbf{E}\left[ \sum_{k=1}^N \tilde X_k(X_k)\right] = \frac{l-1}{1 - \beta} \mathbf{E}\left[\sum_{k=1}^N X_k\right] =  \frac{l-1}{1 - \beta} \alpha.\]

Therefore, 
\begin{align*}
 \mathbf{E}[\min(\sum_{k=1}^N X_k, 1 - \beta)] 
\geq &\mathbf{E}\left[\min\left(\sum_{k=1}^N \tilde X_k(X_k), 1 - \beta\right)\right] \\
= &\mathbf{E}[\min(N'(1 - \beta)/(l-1), 1 - \beta)] \\
= &\frac{1 - \beta}{l-1}\mathbf{E}[\min(N',l-1)]. 
\end{align*}
When $l=2$, this equals $(1 - \beta)\left(1 - e^{-\alpha/ (1 - \beta)}\right)$.

\halmos
\endproof
We are now ready to prove the competitive ratio of $LS$.  The idea is to compare the utilization of each resource $j$ under $LS$ with the utilization of resource $j$ under $OFF$.  The latter is given by $U_j^S + U_j^L$.  The former depends on the choice of the set reserved for $j$, either $L_j$ or $S_j$.  With either choice, we can gauge the total expected utilization, in some cases using Lemma \ref{lm:bound}, to obtain a lower bound.   We then repeat this comparison for all resources $j$ to arrive at a global bound.
\begin{theorem}\label{thm:main}
$LS$ is at least $(1-1/e)/2$-competitive.
\end{theorem}

\noindent \textbf{Proof of Theorem \ref{thm:main}.}
\proof{Proof.}
For each resource $j \in [m]$ there are two cases.
\begin{itemize}
\item Case 1: $U_j^L \geq U_j^S$. Let $Y_j^L$ be the total number of customers who are routed to resource $j$ and whose types are in $L_j$. $Y_j^L$ is a Poisson random variable with mean
\[ \mu_j^L \equiv \mathbf{E}[Y_j^L] = \sum_{i \in L_j} x_{ij}^*.\]
Conditional on the value of $Y_j^L$, the amount of resource $j$ requested by each of the $Y_j^L$ customers is i.i.d. and has mean
\[ \bar u_j^L \equiv \frac{\sum_{i \in L_j} x_{ij}^* u_{ij}}{\sum_{i \in L_j} x_{ij}^*} = \frac{U_j^L}{\mu_j^L}. \]

If $Y_j^L = 1$, the expected amount of resource $j$ taken by that only customer is just $\bar u_j^L$. Thus, we get an expected reward $P(Y_j^L = 1)\bar u_j^L$ from the event $Y_j^L = 1$.

If $Y_j^L>1$, only the first customer can take resource $j$, and all the other $Y_j^L-1$ customers will be rejected due to lack of remaining capacity. The expected amount of resource taken by the first customer may not be $\bar{u}_j^L$ since arrivals are non-homogeneous, but must be still greater than $c_j/2$. Thus, we get an expected reward $P(Y_j^L > 1) c_j/2$ from the event $Y_j^L > 1$.

 In sum, the expected amount of resource taken by these $Y_j^L$ customers is at least
\begin{align}
& P(Y_j^L = 1) \bar u_j^L + P(Y_j^L > 1) c_j/2 \nonumber \\
= & \mu_j^L e^{- \mu_j^L} \bar u_j^L + \left(1 - e^{-\mu_j^L} - \mu_j^L e^{-\mu_j^L}\right) c_j/2 \nonumber \\
= & U_j^L e^{-\mu_j^L}+ \left(1 - e^{-\mu_j^L} - \mu_j^L e^{-\mu_j^L}\right) c_j/2.\label{eq:alg1proofa}
\end{align}

We obtain a lower bound on (\ref{eq:alg1proofa}) by minimizing its value with respect to $\mu_j^L$. We can deduce that
\[ \frac{d}{d\mu_j^L} \left[ U_j^L e^{-\mu_j^L}+ \left(1 - e^{-\mu_j^L} - \mu_j^L e^{-\mu_j^L}\right) c_j/2\right]=0\]
\[\Longrightarrow -U_j^L e^{-\mu_j^L}+ \left(1 + e^{-\mu_j^L}-e^{-\mu_j^L} + \mu_j^L e^{-\mu_j^L}\right) c_j/2 = 0\]
\begin{equation}\label{eq:alg1proofR1a}
 \Longrightarrow \mu_j^L = 2U_j^L  /c_j.
 \end{equation} 

	It is easy to check that (\ref{eq:alg1proofa}) is minimized at solution (\ref{eq:alg1proofR1a}), and the corresponding minimum value of (\ref{eq:alg1proofa}) is
\begin{align*}
& U_j^L e^{-2U_j^L  /c_j}+ \left(1 - e^{-2U_j^L  /c_j} - 2U_j^L  /c_j e^{-2U_j^L  /c_j}\right) c_j/2\\
=& \left(1 - e^{-2U_j^L/c_j}\right) c_j/2\\
\geq & \left( 1 - e^{ - U_j / c_j}\right) c_j/2\\
\geq & \left( 1 - e^{-1}\right) U_j / 2.
\end{align*}
The last step follows since $U_j/c_j \leq 1$.

\item Case 2: $U_j ^L < U_j^S$. Let $Y_j^S$ be the total number of customers who are routed to resource $j$ and whose types are in $S_j$. $Y_j^S$ is a Poisson random variable with mean
\[ \mathbf{E}[Y_j^S] = \sum_{i \in S_j} x_{ij}^*.\]

Let $W^S_1,W^S_2,W^S_3,...,$ be a sequence of i.i.d. random variables each having distribution
\[ P(W^S_k \leq x) = \sum_{i \in S_j} \mathbf{1}(u_{ij} \leq x) \frac{x_{ij}^*}{\sum_{l \in S_j} x_{lj}^*}, \,\,\, \forall k = 1,2,....\]
Here each $W^S_k$ can be seen as the random amount of resource $j$ requested by one of the $Y_j^S$ customers conditional on the value of $Y_j^S$. Then $\sum_{k=1}^{Y_j^S} W^S_k$ represents the total random amount of resource $j$ requested by all the $Y_j^S$ customers. It is easy to check that 
\[ \mathbf{E}\left[\sum_{k=1}^{Y_j^S} W^S_k\right] = U_j^S \geq \frac{U_j}{2}\] 
\[ \Longrightarrow \mathbf{E}\left[\sum_{k=1}^{Y_j^S} \frac{W^S_k}{c_j}\right] = \frac{U_j^S}{c_j} \geq \frac{U_j}{2c_j}.\]

If $\sum_{k=1}^{Y_j^S} W^S_k \leq c_j$, all the $Y_j^S$ customers will be accepted, and we will get total reward $\sum_{k=1}^{Y_j^S} W^S_k$ from resource $j$.
 
  If $\sum_{k=1}^{Y_j^S} W^S_k > c_j$, some of the $Y_j^S$ customers must be rejected due to lack of capacity. But whenever a customer is rejected, the remaining available capacity of resource $j$ must be strictly less than $0.5 c_j$, since $W^S_k/c_j \in [0,0.5]$ w.p.1 for every $k$.
 
 In sum, the total reward we get from resource $j$ is at least 
\begin{align*}
&\sum_{k=1}^{Y_j^S} W^S_k \cdot \mathbf{1}(\sum_{k=1}^{Y_j^S} W^S_k \leq c_j) + 0.5 c_j \cdot \mathbf{1}(\sum_{k=1}^{Y_j^S} W^S_k> c_j)\\
\geq & \min(\sum_{k=1}^{Y_j^S} W^S_k, 0.5 c_j).
\end{align*}
Its expected value can be written as
\[ \mathbf{E}\left[\min\left(\sum_{k=1}^{Y_j^S} W^S_k, \frac{c_j}{2}\right)\right] =  c_j \mathbf{E}\left[\min\left(\sum_{k=1}^{Y_j^S} \frac{W^S_k}{c_j}, \frac{1}{2}\right)\right].\]

 We then apply Lemma \ref{lm:bound} to obtain
\[ \mathbf{E}\left[\min\left(\sum_{k=1}^{Y_j^S} \frac{W^S_k}{c_j}, \frac{1}{2}\right)\right] \geq \frac{1}{2} \left(1 - e^{-2U_j^S/c_j}\right) \geq \frac{1}{2}\left(1-e^{-\frac{U_j}{c_j}}\right)\]
\[ \Longrightarrow \mathbf{E}\left[\min\left(\sum_{k=1}^{Y_j^S} \frac{W^S_k}{c_j}, \frac{1}{2}\right)\right] \geq \frac{U_j}{2c_j}\left(1-\frac{1}{e}\right)\]
\[ \Longrightarrow \mathbf{E}\left[\min\left(\sum_{k=1}^{Y_j^S} W^S_k, \frac{c_j}{2}\right)\right] \geq \frac{U_j}{2}\left(1-\frac{1}{e}\right).\]
\end{itemize}

In sum, in both cases the expected amount of resource $j$ allocated to customers is at least $U_j(1-1/e)/2$. Summing over every resource $j \in [m]$, we can obtain the performance guarantee of our algorithm
\[ \sum_{j \in [m]} U_j \cdot \frac{1}{2} \left(1 - \frac{1}{e}\right) = V^{LP} \cdot \frac{1}{2} \left(1 - \frac{1}{e}\right) \geq \mathbf{E}\left[\OPTI\right]\cdot \frac{1}{2} \left(1 - \frac{1}{e}\right).\]
\halmos
\endproof
Below, we prove that the competitive ratio $(1-1/e)/2$ is tight for $LS$ in the sense that there is at least one problem instance and a corresponding LP solution where $LS$ has relative performance \va{that is bounded above by} $(1-1/e)/2$.
\begin{proposition}
The competitive ratio $(1-1/e)/2$ is tight for $LS$.
\end{proposition}
\proof{Proof.}
We prove the theorem by constructing a special case where the total expected reward of $LS$ is $\mathbf{E}[\OPTI]\cdot (1-1/e)/2$.

Let $n = 2$. Let $\epsilon>0$ be a small number. Let $u_{11} = 0.1$ and $u_{1j} = 0.1(1-\epsilon)$ for every $j\not=1$. Let $u_{21} = 0.5+\epsilon$ and $u_{2j} = (0.5+\epsilon)(1-\epsilon)$ for every $j\not=1$. Let $c_j = 1$ for every $j\in[m]$. The expected number of arrivals is $\Lambda_1 = 5$ and $\Lambda_2 = 0.5/(0.5+\epsilon)$. 

We set $m$ to be very large such that (with probability very close to $1$) the optimal offline algorithm can assign each demand unit to a distinct resource. 
Thus we have
\[ \mathbf{E}[\OPTI] \geq \Lambda_1 u_{12} + \Lambda_2 u_{22} = 5 \cdot 0.1(1-\epsilon) + 0.5/(0.5+\epsilon) \cdot (0.5 + \epsilon)(1-\epsilon) =1-\epsilon.\]

On the other hand, \x{it is easy to verify that the LP has a unique optimal solution, namely  $x_{11}^* = \Lambda_1$, $x_{21}^* = \Lambda_2$ and $x_{ij}^* = 0$ for all other $i, j$.} Then $LS$ reserves resource $1$ only for type-2 customers. In this way, the probability that $LS$ accepts one customer of type 2 is $1 - e^{-\Lambda_2}$. Thus, the expected total reward of $LS$ is
\[ u_{21} \cdot (1 - e^{-\Lambda_2}) = (0.5 + \epsilon) \cdot ( 1- e^{-0.5/(0.5 + \epsilon)}).\]
In sum, the performance ratio of $LS$ can be upper-bounded by
\[ \frac{(0.5 + \epsilon) \cdot ( 1- e^{-0.5/(0.5 + \epsilon)})}{ \mathbf{E}[\OPTI]} \leq \frac{(0.5 + \epsilon) \cdot ( 1- e^{-0.5/(0.5 + \epsilon)})}{1-\epsilon}, \]
which approaches $(1 - 1/e)/2$ when $\epsilon$ tends to $0$. Thus, the competitive ratio of $LS$ is at most $(1 - 1/e)/2$ when $\epsilon$ tends to $0$.
\\
\halmos
\endproof


\x{\subsection{Parameter-dependent bounds for smaller utilization values}}

In this section, we modify the $LS$ algorithm for cases in which all utilization values $u_{ij}$ are bounded away from $c_j$. We characterize the performance of the modified algorithm using competitive ratios that depend on model parameters.

Assume that there is some integer $d \geq 2$ for which $u_{ij} \leq c_j/d$ for all $i \in [n], j \in [m]$. 
The modified $LS$ algorithm considers the value of $u_{ij}$ as ``large'' if $u_{ij} \in (\frac{c_j}{d+1}, \frac{c_j}{d}]$, and ``small'' if $u_{ij} \in [0, \frac{c_j}{d+1}]$. More precisely, we re-define
$$L_j = \{ i \in [n]: u_{ij} > c_j / (d+1)\}$$
and
$$S_j = \{i \in [n]: u_{ij} \leq c_j / (d+1)\}.$$
Also re-define $U_j^L$ and $U_j^S$ based on the modified sets $L_j$ and $S_j$, respectively, in the same way as before. 

In addition to the two options of reserving a resource $j$ for $L_j$ or $S_j$, the modified $LS$ algorithm considers a third option of simply pooling all customer types together. In other words, the third option allows for any customer who is routed to a resource $j$ to take the resource in a first-come first-served fashion. 

\x{Intuitively, if we only considered the same (first two) options as in the original $LS$ algorithm, the competitive ratio would never exceed $0.5$, because in the worst case, the algorithm would always reject half of total demand by the reservation rule. By adding the third option, which potentially utilizes all the demand, we can raise the competitive ratio of the modified $LS$ algorithm  to $1$ as $d \to \infty$. Note that, however, this third option is not helpful in the general case $d=1$, when it can cause arbitrarily poor performance to pack all customer types together.}

\noindent {\bf Modified Large-or-Small ($MLS$) Algorithm:}
\begin{enumerate}
\item (Pre-processing step) Same as for $LS$.
\item (Reservation step) Calculate the following three ratios
\begin{align*}
 \ratioL & \equiv \frac{1}{d+1} \left[ \sum_{k=1}^d e^{-\frac{(d+1)U_j^L}{c_j}  } \frac{\left(\frac{(d+1)U_j^L}{c_j}\right)^k}{(k-1)!} + \sum_{k=d+1}^\infty d\cdot e^{-\frac{(d+1)U_j^L}{c_j}} \frac{\left( \frac{(d+1)U_j^L}{c_j} \right)^k}{k!} \right],\\
 \ratioS & \equiv \frac{1}{d+1} \left[ \sum_{k=1}^d e^{-\frac{(d+1)U_j^S}{c_j}} \frac{\left( \frac{(d+1)U_j^S}{c_j}\right)^k}{(k-1)!} + \sum_{k=d+1}^\infty d \cdot e^{-\frac{(d+1)U_j^S}{c_j}} \frac{\left( \frac{(d+1)U_j^S}{c_j}\right)^k}{k!}\right],\\
\ratioall & \equiv \frac{U_j}{c_j} \cdot \left( 1 - e^{-d} \sum_{i=d}^\infty (i-d+1) \frac{d^{i-1}}{i!}\right).
\end{align*}
\x{Reserve resource $j$ for customer types in $L_j$ if $\ratioL$ is the largest among the three ratios; reserve resource $j$ for customer types in $S_j$ if $\ratioS$ is the largest among the three; otherwise, open the resource to all customer types.}
\item (Random routing step) Same as for $LS$.
\item (Admission step) Same as for $LS$.
\end{enumerate}

\begin{figure}[h!]
  \centering
  \caption{Competitive ratio $\max(\ratioL, \ratioS, \ratioall)$ of $MLS$ as a function of $U_j^S \in [0, U_j]$ when $U_j = c_j = 1$. In particular, $\ratioall$ is independent of $U_j^S$ as we fix $U_j=1$. $\ratioL$ depends on $U_j^S$ as $U_j^L = U_j - U_j^S$. The flat part in the middle of each curve corresponds to $\ratioall$ being larger than $\ratioL$ and $\ratioS$. The decreasing part of each curve corresponds to $\ratioL$ being the largest (when $U_j^S$ is small and $U_j^L$ is large).  The increasing part corresponds to $\ratioS$ being the largest.}
  \label{fig:ratios}
  \includegraphics[width=0.5\textwidth]{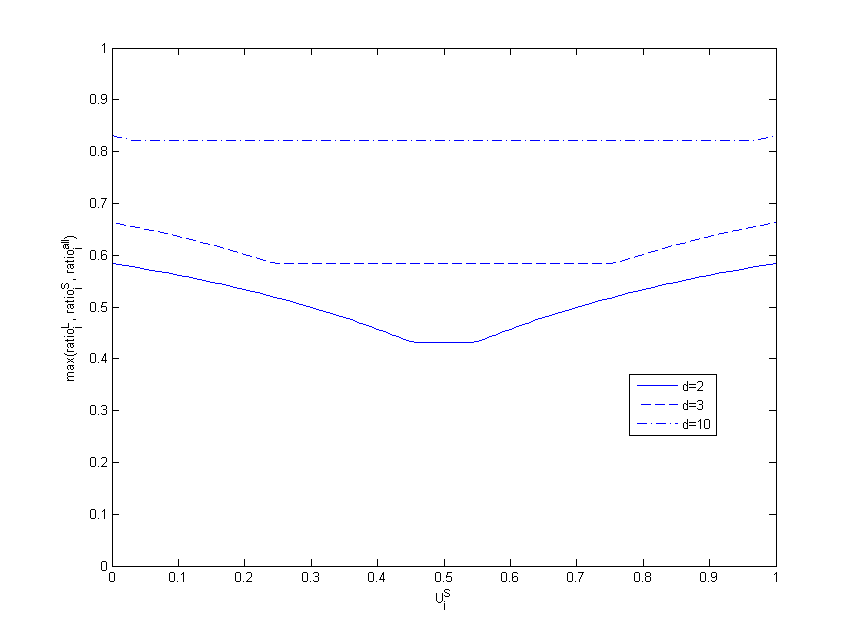}
\end{figure}

\x{Theorem \ref{thm:mls} states that the performance of the $MLS$ algorithm is determined by $\max\{\ratioL, \ratioS, \ratioall\}$, which is essentially a function of $d$, $U_j^L$ and $U_j^S$. Each of the three ratios represents the competitive ratio when the algorithm uses the corresponding reservation strategy. Figure \ref{fig:ratios} illustrates when each of the three ratios is the largest.
For example, when $d=2$, $U_j = c_j$ and $U_j^S = 0.2 U_j$, we have $\ratioS \approx 0.18$, $\ratioL \approx 0.54$ and $\ratioall \approx 0.43$, while the upper bound given by Proposition \ref{prop:upperboundd} is about $0.86$.}

\x{We have $\ratioL$ increases in $U_j^L$, and $\ratioS$ increases in $U_j^S$, because intuitively, the reservation strategies perform better when there are more customers reserved. When $U_j^L$ and $U_j^S$ are balanced, $\ratioall$ is the largest (for $d \geq 2$), because this option pools all customers together and does not discard half of total demand.  }

\begin{theorem}\label{thm:mls}
 Suppose there is some integer $d\geq 2$ for which $u_{ij}\leq \frac{c_j}{d}$ for all $i$ and $j$. From each resource $j$, the modified $LS$ algorithm earns expected reward
\[ c_j \max(\ratioL, \ratioS, \ratioall).\]
\end{theorem}

\proof{Proof.}
\x{Based on the modified sets $S_j$ and $L_j$, re-define variables such as $\yjS$, $\yjL$, $\mujL$, $\bar u_j^L$ and $\{W_k^S\}_{k \geq 1}$ in the same way as in the proof of Theorem \ref{thm:main}.}

\begin{itemize}
\item Case 1: Suppose the modified algorithm reserves resource $j$ for customer types in $L_j$, i.e., $\ratioL = \max(\ratioL, \ratioS, \ratioall)$ .

\x{The algorithm will allocate resource $j$ to exactly $\min(d, \yjL)$ customers whose types are in $L_j$. Conditioned on $\yjL \leq d$, each of these $\min(d, \yjL) = \yjL$ customers will take $\bar u_j^L = U_j^L / \mujL$ unit of resource $j$ in expectation. Conditioned on $\yjL > d$, each of these $\min(d, \yjL)=d$ customers must take at least $\frac{c_j}{d+1}$ unit of resource $j$. Thus, We expect to earn at least}
\begin{align}
& \sum_{k=1}^d k\bar u_j^L \cdot P(Y_j^L = k)  + \frac{d c_j}{d+1} \cdot P(Y_j^L > d) \label{eq:r2thm2pfcase1b} \\
= & \sum_{k=1}^d k \frac{U_j^L}{\mujL} \cdot e^{-\mujL} \frac{(\mujL)^k}{k!} + \sum_{k=d+1}^\infty \frac{dc_j}{d+1} \cdot e^{-\mujL} \frac{(\mujL)^k}{k!}\nonumber \\
= & \sum_{k=1}^d  U_j^L \cdot e^{-\mujL} \frac{(\mujL)^{k-1}}{(k-1)!} + \sum_{k=d+1}^\infty \frac{dc_j}{d+1} \cdot e^{-\mujL} \frac{(\mujL)^k}{k!}\nonumber\\ 
= & U_j^L P(\yjL < d) + \frac{dc_j}{d+1} P(\yjL > d). \nonumber
 \end{align}
 
 We want to find the $\mujL = \bE[\yjL]$ that minimizes this expected reward. We can deduce that
 \begin{align*}
&  \frac{\partial}{\partial \mujL} [U_j^L P(\yjL < d) + \frac{dc_j}{d+1} P(\yjL > d)] \\
  =& U_j^L P(\yjL < d-1) - U_j^L P(\yjL < d) + \frac{dc_j}{d+1} P(\yjL > d-1)  - \frac{dc_j}{d+1} P(\yjL > d)\\
  =& \frac{dc_j}{d+1} P(\yjL = d) - U_j^L P(\yjL = d-1)\\
  =& \frac{dc_j}{d+1} e^{-\mujL} \frac{(\mujL)^d}{d!} - U_j^L e^{-\mujL}\frac{(\mujL)^{d-1}}{(d-1)!}\\
  =& e^{-\mujL}\frac{(\mujL)^{d-1}}{(d-1)!} \left( \frac{dc_j}{d+1} \cdot \frac{\mujL}{d}  - U_j^L\right).
  \end{align*}
Setting the derivative to zero, we obtain $\mujL = (d+1) \frac{U_j^L}{c_j}$. It is easy to check that \eqref{eq:r2thm2pfcase1b} is minimized at $\mujL= (d+1) \frac{U_j^L}{c_j}$, when it is equal to
\begin{align*}
&\sum_{k=1}^d k \frac{U_j^L}{\mujL} \cdot e^{-\mujL} \frac{(\mujL)^k}{k!} + \sum_{k=d+1}^\infty \frac{dc_j}{d+1} \cdot e^{-\mujL} \frac{(\mujL)^k}{k!}\\
= & \sum_{k=1}^d \frac{k c_j}{d+1} \cdot e^{-\frac{(d+1)U_j^L}{c_j}  } \frac{\left(\frac{(d+1)U_j^L}{c_j}\right)^k}{k!} + \sum_{k=d+1}^\infty \frac{dc_j}{d+1} \cdot e^{-\frac{(d+1)U_j^L}{c_j}} \frac{\left( \frac{(d+1)U_j^L}{c_j} \right)^k}{k!}\\
= & \frac{c_j}{d+1} \left[ \sum_{k=1}^d e^{-\frac{(d+1)U_j^L}{c_j}  } \frac{\left(\frac{(d+1)U_j^L}{c_j}\right)^k}{(k-1)!} + \sum_{k=d+1}^\infty d\cdot e^{-\frac{(d+1)U_j^L}{c_j}} \frac{\left( \frac{(d+1)U_j^L}{c_j} \right)^k}{k!} \right]\\
= & c_j \ratioL.
\end{align*}

Therefore, when the algorithm reserves resource $j$ for $L_j$, $c_j \ratioL = c_j \max(\ratioL, \ratioS, \ratioall)$ lower bounds the expected reward from resource $j$.

\item Case 2: Next, suppose the algorithm reserves resource $j$ for customer types in $S_j$, i.e., $\ratioS = \max(\ratioL, \ratioS, \ratioall)$.

\x{Recall that $\{\wkS\}_{k\geq 1}$ is a sequence of i.i.d. random variables, each representing the random amount of resource $j$ requested by one of the $\yjS$ customers, conditional on the value of $\yjS$.} Since $u_{ij} \leq c_j / (d+1)$ for $i \in S_j$, we must have
\[ \wkS / c_j \in [0, \frac{1}{d+1}].\]
Conditioned on $\yjS$, if $\sum_{k=1}^\yjS \wkS > c_j$, the remaining capacity of resource $j$ must be strictly less than $\frac{c_j}{d+1}$ by the end. Thus the total reward we  expect to earn from resource $j$ is
\begin{align*}
& \bE\left[ \min\left( \sum_{k=1}^\yjS \wkS, c_j - \frac{c_j}{d+1} \right) \right]\\
= & c_j \cdot \bE\left[ \min\left( \sum_{k=1}^\yjS \frac{\wkS}{c_j}, 1 - \frac{1}{d+1} \right) \right]\\
\geq & c_j \cdot \frac{1}{d+1} \left[ \sum_{k=1}^d e^{-\frac{(d+1)U_j^S}{c_j}} \frac{\left( \frac{(d+1)U_j^S}{c_j}\right)^k}{(k-1)!} + \sum_{k=d+1}^\infty d \cdot e^{-\frac{(d+1)U_j^S}{c_j}} \frac{\left( \frac{(d+1)U_j^S}{c_j}\right)^k}{k!}\right]\\
= & c_j \ratioS.
\end{align*}
Here the inequality is by Lemma \ref{lm:bound}, when $\beta = 1/(d+1)$, $l = d+1$ and $\alpha = U_j^S / c_j$. Therefore, when the algorithm reserves resource $j$ for $S_j$, $c_j \ratioS = c_j \max(\ratioL, \ratioS, \ratioall)$ lower bounds the expected reward from resource $j$.

\item Case 3: Finally, suppose the algorithm opens resource $j$ to all customer types. 
Let $Y_j$ be the total number of customers who are routed to resource $j$. Then $Y_j$ is a Poisson random variable with mean
\[ \mathbf{E}[Y_j] = \sum_{i\in [n]} x_{ij}^*.\]

Let $W_1,W_2,W_3,...,$ be a sequence of i.i.d. random variables each having distribution
\[ P(W_k \leq x) = \sum_{i\in [n]} \mathbf{1}(u_{ij} \leq x) \frac{x_{ij}^*}{\sum_{l \in [n]} x_{lj}^*}, \,\,\, \forall k = 1,2,....\]
Here each $W_k$ can be seen as the random amount of resource $j$ requested by one of the $Y_j$ customers conditional on the value of $Y_j$. Then $\sum_{k=1}^{Y_j} W_k$ represents the total random amount of resource $j$ requested by all the $Y_j$ customers, conditioned on $Y_j$. It is easy to check that $ \mathbf{E}\left[\sum_{k=1}^{Y_j} \frac{W_k}{c_j}\right] = \frac{U_j}{c_j}$.
 
Since $W_k/c_j \in [0,1/d]$ for every $k$, the expected amount of resource $j$ taken by these $Y_j$ customers is at least 
\[ \mathbf{E}\left[\min\left(\sum_{k=1}^{Y_j} W_k, c_j - c_j/d\right)\right] =  c_j \mathbf{E}\left[\min\left(\sum_{k=1}^{Y_j} \frac{W_k}{c_j}, 1 - 1/d\right)\right].\]

We then apply Lemma \ref{lm:bound} to obtain

\[ c_j \mathbf{E}\left[\min\left(\sum_{k=1}^{Y_j} \frac{W_k}{c_j}, 1 - 1 /d\right)\right] \geq c_j \frac{1 - 1/d}{d-1} \mathbf{E}[\min(N',d-1)] = \frac{c_j}{d} \mathbf{E}[\min(N',d-1)]  ,\]
where $N'$ is a Poisson random variable with mean $\frac{U_j}{c_j}\cdot \frac{d-1}{1-1/d} = \frac{U_j}{c_j} \cdot d$. Let $N$ be a Poisson random variable with mean $\mathbf{E}[N] = d$ that is independent of $N'$. We have $\mathbf{E}[N] \geq \mathbf{E}[N']$ since $U_j / c_j \leq 1$.
We can further deduce that

\begin{align*}
& \frac{c_j}{d} \mathbf{E}[\min(N',d-1)]\\
= & \frac{c_j}{d} \left[ \sum_{i=1}^{d-1} i \cdot e^{-\mathbf{E}[N']} \frac{(\mathbf{E}[N'])^i}{i!} + \sum_{i=d}^\infty (d-1) e^{-\mathbf{E}[N']} \frac{(\mathbf{E}[N'])^i}{i!}\right]\\
= & \frac{c_j}{d} \mathbf{E}[N'] \left[ \sum_{i=0}^{d-2} e^{-\mathbf{E}[N']} \frac{(\mathbf{E}[N'])^i}{i!} + \sum_{i=d-1}^\infty \frac{d-1}{i+1} \cdot e^{-\mathbf{E}[N']} \frac{(\mathbf{E}[N'])^i}{i!}\right]\\
= & \frac{c_j}{d} \mathbf{E}[N'] \mathbf{E}[\min(1,\frac{d-1}{N'+1})]\\
\geq & \frac{c_j}{d} \mathbf{E}[N'] \mathbf{E}[\min(1,\frac{d-1}{N+1})]\\
= & \frac{c_j}{d} \frac{\mathbf{E}[N']}{\mathbf{E}[N]} \cdot \mathbf{E}[N]\mathbf{E}[\min(1,\frac{d-1}{N+1})]\\
= & \frac{c_j}{d} \frac{\mathbf{E}[N']}{\mathbf{E}[N]} \cdot  \mathbf{E}[N] \left[ \sum_{i=0}^{d-2} e^{-\mathbf{E}[N]} \frac{(\mathbf{E}[N])^i}{i!} + \sum_{i=d-1}^\infty \frac{d-1}{i+1} \cdot e^{-\mathbf{E}[N]} \frac{(\mathbf{E}[N])^i}{i!}\right]\\
= & \frac{c_j}{d} \frac{\mathbf{E}[N']}{\mathbf{E}[N]} \cdot \left[ \sum_{i=1}^{d-1} i \cdot e^{-\mathbf{E}[N]} \frac{(\mathbf{E}[N])^i}{i!} + \sum_{i=d}^\infty (d-1) e^{-\mathbf{E}[N]} \frac{(\mathbf{E}[N])^i}{i!}\right]\\
= & \frac{c_j}{d} \frac{\mathbf{E}[N']}{\mathbf{E}[N]} \mathbf{E}[\min(N, d-1)] \\
= & U_j \cdot \frac{1}{d}  \mathbf{E}[\min(N,d-1)]\\
= & U_j \cdot \frac{1}{d} (\mathbf{E}[N] - \mathbf{E}[\max(N-d+1,0)])\\
= & U_j \cdot \left( 1 - e^{-d} \sum_{i=d}^\infty (i-d+1) \frac{d^{i-1}}{i!}\right)\\
= &  c_j \ratioall.
\end{align*}
Therefore, when the algorithm reserves resource $j$ for all customer types, $c_j \ratioall = c_j \max(\ratioL, \ratioS, \ratioall)$ lower bounds the expected reward from resource $j$.
\end{itemize}
\halmos
\endproof


\x{As illustrated in Figure \ref{fig:ratios}, $\ratioall$ increases in $d$. In the next corollary, we prove that $\ratioall = 1 - O(1/\sqrt{d})$. This matches the best possible dependence on $d$ according to Proposition \ref{prop:upperboundd2}.}

\x{\begin{corollary} \label{cor:improveBounded} If there is some integer $d\geq 2$ for which $u_{ij}\leq \frac{c_j}{d}$ for all $i$ and $j$, then the competitive ratio of the modified $LS$ is at least $1 - \frac{1}{\sqrt{2\pi d}} + O(1/d)$.
\end{corollary}}
\proof{Proof.}  
The expected amount of resource $j$ allocated to customers is at least
\begin{align*}
& c_j \max(\ratioL, \ratioS, \ratioall)\\
\geq & c_j \ratioall\\
= & U_j \cdot \left( 1 - e^{-d} \sum_{i=d}^\infty (i-d+1) \frac{d^{i-1}}{i!}\right)\\
= & U_j \cdot \left( 1 - e^{-d} \frac{d^d}{d!} - \frac{1}{d} e^{-d} \sum_{i=d}^\infty \frac{d^i}{i!} \right)\\
\geq & U_j \cdot \left( 1 - e^{-d} \frac{d^d}{d!} - \frac{1}{d} e^{-d} \sum_{i=0}^\infty \frac{d^i}{i!} \right)\\
= & U_j \cdot \left( 1 - e^{-d} \frac{d^d}{d!} - \frac{1}{d}\right)\\
= & U_j \cdot \left( 1 - \frac{1}{\sqrt{2\pi d}} + O(1/d)\right) .
\end{align*}
The last step above follows by Stirling's formula.

Summing over every resource $j \in [m]$, we can obtain the performance guarantee of the modified $LS$ algorithm
\begin{align*}
 &\sum_{j \in [m]} U_j \cdot \left( 1 - \frac{1}{\sqrt{2\pi d}} + O(1/d)\right) \\
 = & V^{LP} \cdot \left( 1 - \frac{1}{\sqrt{2\pi d}} + O(1/d)\right)\\
 \geq & \mathbf{E}\left[\OPTI\right]\cdot \left( 1 - \frac{1}{\sqrt{2\pi d}} + O(1/d)\right).
 \end{align*}
\halmos
\endproof


\section{Improving the Bound}\label{sec:improve}
In this section, we derive an algorithm with an improved competitive ratio compared to $LS$. This algorithm also groups customer types based on the utilization $u_{ij}$, but in a more sophisticated way than $LS$. Moreover, this algorithm also relaxes the random routing step in order to allow customers more opportunities to be assigned to resources. This strategy of allowing greater resource sharing among customer types greatly improves the empirical performance of the algorithm.

We will prove that the competitive ratio of the new algorithm is
\begin{equation}
\label{eq:cra}
 r^* =  \max\left\{ r \in (0,0.5) :  r \leq \max_{z \in (0,0.5)} h(z,r) \right\},
 \end{equation}
 where 
\begin{equation}
 h(z,r) \equiv  z - \left[ z - \frac{1}{2} \left( 1 - \frac{1}{1-2r}\cdot \frac{1}{e^2} \right) \right] (1 - 2r) \left( \frac{1 - z}{1 - z - r}\right) ^{2(1-z)}. \end{equation}
We can numerically solve (\ref{eq:cra}) to find that $r^* \approx 0.321$.

For every resource $j$, we first divide all customer types into two sets $S_j$ and $L_j$ in the same way that $LS$ does. Then, we further partition the customers in $S_j$ into two sets, ``medium small" and ``tiny", depending on their utilization of resource $j$. Let
\begin{eqnarray} M_j &=& \{ i \in S_j : u_{ij} \geq z^* \cdot c_j\} ,\\
T_j &=& \{ i \in S_j : u_{ij} < z^* \cdot c_j \},\end{eqnarray}
where
 \begin{equation} z^* \equiv \arg\!\!\!\max_{z \in (0,0.5)}  h(z,r^*) \approx 0.42. \end{equation}
It is easy to check that there is only one maximizer.

Recall that $U_j = \sum_{i \in [n]} x_{ij}^*u_{ij}$, $U_j^L = \sum_{i \in L_j} x_{ij}^*u_{ij}$ and $U_j^S = \sum_{i \in S_j} x_{ij}^* u_{ij}$, where $x^*$ is an optimal solution to the LP (\ref{eq:lp}). We further define $U_j^M \equiv \sum_{i \in M_j} x_{ij}^*u_{ij}$ and $U_j^T \equiv \sum_{i \in T_j} x_{ij}^*u_{ij}$ analogously.

Intuitively, the load values $U_j^T, U_j^M, U_j^S, U_j^L$ serve as estimates for how much capacity of resource $j$ is expected to be utilized by customers of types in the sets $T_j, M_j, S_j$ and $L_j$, respectively. For a given resource $j$, if any of the load values dominates the others, it might be a good strategy to reserve resource $j$ exclusively for customers in the corresponding set. 

We next categorize every resource $j$ into one of two types based on the load values.
\begin{definition} \label{def:resourcetype}
Resource $j$ is a type-A resource if 
\[ U_j^S \geq -  0.5 c_j  \log(1 - 2r^*U_j/c_j) \]
\[\mbox{or}\quad U_j^T \geq -  (1 - z^*) c_j \log\left(1 - \frac{r^*U_j}{c_j(1-z^*)}\right). \]
Otherwise, resource $j$ is a type-B resource.
\end{definition}
The motivation for the above definition is as follows.  If resource
$j$ is of type A, then $U_j^S$ or $U_j^T$ are relatively large
compared to other load values, which implies that customers that are
routed to resource $j$ by the LP (\ref{eq:lp}) tend to have
relatively small utilization $u_{ij}$. On the other hand, if resource
$j$ is of type B, then $U_j^L$ is relatively larger, which implies
that customers that are routed to resource $j$ by the LP tend to have
relatively large utilization $u_{ij}$.

Depending on the type of resource, we will reserve each resource wholly for a certain set of customer types. We say that a customer of type $i$ is \emph{admissible} to resource $j$ if this customer can be assigned to resource $j$ by our algorithm.
The following definition defines the reservation criteria of the algorithm.
\begin{definition}
A customer of type $i$, $i \in [n]$, is admissible to resource $j$, $j \in [m]$, if and only if at least one of the following criteria holds:
\begin{itemize}
\item Resource $j$ is of type A,
\item or $i \in M_j \cup L_j$.
\end{itemize}
\end{definition}

We are now ready to specify the improved algorithm.\\
\noindent {\bf Refined Large-or-Small Algorithm ($RLS$):}
\begin{enumerate}
\item (Pre-processing step) Same as for $LS$.
\item (Random routing step) Same as for $LS$.
\item (Admission step) If a customer is admissible to resource $j$ and there is enough remaining capacity in $j$, then assign the customer to resource $j$.    
\item (Resource-sharing step) If a customer is rejected in the Admission Step, but there is another resource with enough remaining capacity and to which the customer is admissible, then assign the customer to any such resource.  Otherwise, reject the customer.
\end{enumerate}

The idea of the algorithm is as follows.  If a resource $j$ is of type
$A$, then we can bound from above the left-over capacity of the
resource, because a type-$A$ resource tends to be used by a
sufficiently high number of customers in sets $S_j$ and
$T_j$. Whenever one such customer is rejected, we know that the
remaining capacity is small. Furthermore, since it is not
disadvantageous to turn away customers of other types, their
utilization being higher than those in $S_j$ and $T_j$, we will admit
customers of all types to a type-$A$ resource.  On the other hand, if
resource $j$ is of type $B$, then customers who are admissible to the
resource have large utilization values. We can allocate a large enough
amount of the resource as soon as one such customer arrives. We do not
admit customers with small utilization values to type-$B$ resources in
order to leave enough space for customer types in $M_j$ and $L_j$.

For each resource $j$, let $N_j$ denote the total number of customers who are routed to resource $j$ in Step (2) of $RLS$. Note that $N_j$ does not include customers who are assigned resource $j$ in Step (4) of $RLS$. Let $W_{j1}, W_{j2},...,$ be a sequence of i.i.d. random variables each having a distribution that is given by
\[ P(W_{j1} \leq x) = \sum_{i \in [n]} \mathbf{1}(u_{ij} \leq x) \frac{x_{ij}^*}{\sum_{k \in [n]} x_{kj}^*}.\]
That is, each variable $W_{j1}$ can be seen as the utilization $u_{ij}$ of a single random customer who is routed to resource $j$ during the horizon.  Since the probability that such a customer has type $i$ is $\frac{x_{ij}^*}{\sum_{k \in [n]} x_{kj}^*}$, $W_{j1}$ takes value $u_{ij}$ with this probability.  The following lemma gives a lower bound on the expected amount of capacity of a type-A resource that will be allocated by $RLS$.

\begin{lemma}\label{lem:typeAbound}
If resource $j$ is of type-A, then the expected amount of resource $j$ allocated by $RLS$ is at least
\[ \max\left\{ \mathbf{E}\left[\min( \sum_{k=1}^{N_j} W_{jk}\mathbf{1}( W_{jk} \leq 0.5 c_j), 0.5c_j)\right], \mathbf{E}\left[\min( \sum_{k=1}^{N_j} W_{jk}\mathbf{1}( W_{jk} \leq z^*c_j), (1-z^*)c_j)\right] \right\}   .\]
\end{lemma}
\proof{Proof.}
$\sum_{k=1}^{N_j} W_{jk}\mathbf{1}( W_{jk} \leq 0.5 c_j)$ has the same distribution as the total (random) amount of resource $j$ requested by customers who are routed to resource $j$ and whose types are in $S_j$. If the actual amount of resource $j$ allocated to customers is less than $\sum_{k=1}^{N_j} W_{jk}\mathbf{1}( W_{jk} \leq 0.5 c_j)$, it must be that at least one customer with type in $S_j$ is rejected due to lack of remaining capacity of resource $j$. In such a case, the actual  amount of resource $j$ allocated to customers must be at least $0.5c_j$, because otherwise the customer with type in $S_j$ would not have been rejected. Thus, the total amount of resource $j$ allocated by our algorithm is at least $\min( \sum_{k=1}^{N_j} W_{jk}\mathbf{1}( W_{jk} \leq 0.5 c_j), 0.5c_j)$.

A similar argument applies to customers with types in $T_j$. $\sum_{k=1}^{N_j} W_{jk}\mathbf{1}( W_{jk} \leq z^*c_j)$ has the same distribution as the total amount of resource $j$ requested by customers who are routed to resource $j$ and whose types are in $T_j$. If at least one of these requests is rejected, the remaining capacity of resource $j$ must be at most $z^* c_j$. Thus, the total amount of resource $j$ allocated to customers is at least $\min( \sum_{k=1}^{N_j} W_{jk}\mathbf{1}( W_{jk} \leq z^*c_j), (1-z^*)c_j)$.

The proof follows when we take expectation of the lower bounds.
\halmos
\endproof

Let $\mu_j^M = \sum_{i \in M_j} x_{ij}^*$ and $\mu_j^L = \sum_{i \in L_j} x_{ij}^*$ be the expected number of customers who are routed to resource $j$ and whose types are in $M_j$ and $L_j$, respectively. The following lemma gives a lower bound on the expected amount of capacity of a type-B resource that will be allocated by $RLS$.

\begin{lemma}\label{lem:typeBbound}
If resource $j$ is of type B, then the expected amount of resource $j$ allocated to customers in $RLS$ is at least
\[ \min\{ z^* c_j,  e^{-\mu_j^M}[ U_j^L e^{-\mu_j^L} + 0.5 c_j ( 1 - e^{-\mu_j^L} - \mu_j^L e^{-\mu_j^L})] + (1 - e^{-\mu_j^M}) z^*c_j \}.\]
\end{lemma}
\proof{Proof.}
If any customer is assigned to resource $j$ in Step (4) of the algorithm, then at least $z^* c_j$ of resource $j$ is allocated, since every customer type $i$ that is admissible to resource $j$ satisfies $u_{ij} \geq z^* c_j$. 

If no customer is assigned to resource $j$ by Step (4), then resource $j$ can only be allocated to customers who are directly routed to this resource, i.e. in Step (3) of $RLS$. We consider three cases:
\begin{itemize}
\item At least one customer with type in $M_j$ is routed to resource $j$. This event occurs with probability $1 - e^{-\mu_j^M}$. In such a case, we use $z^*c_j$ as the lower bound on the amount of resource $j$ taken by customers.
\item No customer with type in $M_j$ is routed to resource $j$, and exactly one customer with type in $L_j$ is routed to resource $j$. This event occurs with probability $e^{-\mu_j^M}  \cdot \mu_j^L e^{-\mu_j^L}$. Conditional on this event, the expected amount of resource $j$ taken by the only customer with type in $L_j$ is
\[ \frac{\sum_{i \in L_j} x_{ij}^* u_{ij}}{\sum_{i \in L_j} x_{ij}^*} = \frac{U_j^L}{\mu_j^L}.\]
\item No customer with type in $M_j$ is routed to resource $j$, and more than one customer with type in $L_j$ are routed to resource $j$. This event occurs with probability $e^{-\mu_j^M}( 1 - e^{-\mu_j^L} - \mu_j^L e^{-\mu_j^L})$. In this event, we use $0.5c_j$ as the lower bound on the amount of resource $j$ taken by the customer in $L_j$, by definition of $L_j$. 
\end{itemize}
In summary, if no customer is assigned to resource $j$ in Step (3), the expected amount of resource $j$ taken by customers routed to the resource is at least
\begin{align*}
& (1 - e^{-\mu_j^M}) z^* c_j + e^{-\mu_j^M}  \cdot \mu_j^L e^{-\mu_j^L} \cdot \frac{U_j^L}{\mu_j^L} + e^{-\mu_j^M}( 1 - e^{-\mu_j^L} - \mu_j^L e^{-\mu_j^L}) 0.5 c_j \\
= & e^{-\mu_j^M}[ U_j^L e^{-\mu_j^L} + 0.5 c_j ( 1 - e^{-\mu_j^L} - \mu_j^L e^{-\mu_j^L})] + (1 - e^{-\mu_j^M}) z^*c_j.
\end{align*}

We complete the proof by combining this result with the lower bound $z^* c_j$ for the case that a customer is assigned to resource $j$ by Step (4) of the algorithm.
\halmos
\endproof

We combine the previous two lemmas to prove the performance guarantee of the algorithm.
\begin{theorem}
For each resource $j$, the expected amount of resource $j$ allocated to customers is at least $r^* U_j$.
\end{theorem}
\proof{Proof.}
First consider the case that resource $j$ is of type A. Since $\sum_{k=1}^{N_j} W_{jk}\mathbf{1}( W_{jk} \leq 0.5 c_j)$ and $\sum_{k=1}^{N_j} W_{jk}\mathbf{1}( W_{jk} \leq z^*c_j)$ has compound Poisson distribution with mean
\[ \mathbf{E}\left[\sum_{k=1}^{N_j} W_{jk}\mathbf{1}( W_{jk} \leq 0.5 c_j)\right] = U_j^S, \,\,\,\,\,\,\, \mathbf{E}\left[\sum_{k=1}^{N_j} W_{jk}\mathbf{1}( W_{jk} \leq z^*c_j)\right] = U_j^T,\]
we can apply Lemma \ref{lm:bound} to get
\begin{align*}
 \mathbf{E}\left[ \min( \sum_{k=1}^{N_j} W_{jk}\mathbf{1}( W_{jk} \leq 0.5 c_j), 0.5c_j)\right] 
= & c_j \mathbf{E}\left[ \min( \sum_{k=1}^{N_j} \frac{W_{jk}}{c_j} \mathbf{1}( W_{jk} \leq 0.5c_j), 0.5)\right] \\
\geq & c_j 0.5 \left(1 - e^{-2U_j^S/c_j}\right),
\end{align*}
and
\begin{align*}
 \mathbf{E}\left[ \min( \sum_{k=1}^{N_j} W_{jk}\mathbf{1}( W_{jk} \leq z^*c_j), (1-z^*)c_j)\right] 
= & c_j \mathbf{E}\left[ \min( \sum_{k=1}^{N_j} \frac{W_{jk}}{c_j} \mathbf{1}( W_{jk} \leq z^*c_j), 1-z^*)\right] \\
\geq & c_j (1-z^*) \left(1 - e^{-\frac{U_j^T}{c_j(1-z^*)}}\right).
\end{align*}

Then according to Definition \ref{def:resourcetype}, we have by the definition of type-A resources
\[ U_j^S \geq -  0.5 c_j  \log(1 - 2r^*U_j/c_j) \quad \Longrightarrow \quad c_j 0.5 (1 - e^{-2U_j^S/c_j}) \geq r^* U_j,\]
or
\[ U_j^T \geq -  (1 - z^*) c_j \log(1 - \frac{r^*U_j}{c_j(1-z^*)})\quad \Longrightarrow \quad c_j (1-z^*) (1 - e^{-\frac{U_j^T}{c_j(1-z^*)}}) \geq r^*U_j.\]

In sum, we have 
\begin{align*}
& \max \left\{\mathbf{E}\left[ \min( \sum_{k=1}^{N_j} W_{jk}\mathbf{1}( W_{jk} \leq 0.5 c_j), 0.5c_j)\right] , \mathbf{E}\left[ \min( \sum_{k=1}^{N_j} W_{jk}\mathbf{1}( W_{jk} \leq z^*c_j), (1-z^*)c_j)\right] \right\} \\
\geq & \max \left\{ c_j 0.5 (1 - e^{-2U_j^S/c_j}), c_j (1-z^*) (1 - e^{-\frac{U_j^T}{c_j(1-z^*)}}) \right\}\\
\geq & r^*U_j. 
\end{align*}
This proves the theorem for type-A resources.

Now we consider the case that resource $j$ is of type B. Starting from this point, we will assume without loss of generality that $c_j = 1$. Based on Lemma \ref{lem:typeBbound}, we need to show
\[ \min\{ z^* ,  e^{-\mu_j^M}[ U_j^L e^{-\mu_j^L} + 0.5 ( 1 - e^{-\mu_j^L} - \mu_j^L e^{-\mu_j^L})] + (1 - e^{-\mu_j^M}) z^* \} \geq r^* U_j.\]
Since $z^* > r^*$ as we numerically checked, it suffices to show
\[ e^{-\mu_j^M}\left[ U_j^L e^{-\mu_j^L} + 0.5 \left( 1 - e^{-\mu_j^L} - \mu_j^L e^{-\mu_j^L}\right)\right] + \left(1 - e^{-\mu_j^M}\right) z^* \geq r^* U_j\]
based on Lemma \ref{lem:typeBbound}.

By examining the first and second derivatives of $U_j^L e^{-\mu_j^L} + 0.5  ( 1 - e^{-\mu_j^L} - \mu_j^L e^{-\mu_j^L})$ with respect to $\mu_j^L$, it is easy to check that 
\begin{align*}
e^{-\mu_j^M} \left[ U_j^L e^{-\mu_j^L} + 0.5 \left(1 - e^{-\mu_j^L} - \mu_j^L e^{-\mu_j^L}\right)\right] + \left(1 - e^{-\mu_j^M}\right) z^*
&\geq & e^{-\mu_j^M} \cdot  0.5 \left(1 - e^{-2U_j^L}\right) + \left(1 - e^{-\mu_j^M}\right) z^*\\
&= & z^* - e^{-\mu_j^M}\left[ z^* - 0.5\left(1 - e^{-2U_j^L}\right)\right].
\end{align*}

If $z^* < 0.5(1 - e^{-2U_j^L})$, we must have $z^* - e^{-\mu_j^M}[ z^* - 0.5(1 - e^{-2U_j^L})] > z^* > r^* = r^* c_j  \geq r^* U_j$, which proves the theorem for this case.

Now suppose $z^* \geq 0.5(1 - e^{-2U_j^L})$. Since $\mu_j^M = \sum_{i \in M_j} x_{ij}^* \geq \sum_{i \in M_j} x_{ij}^* \cdot 2 u_{ij} = 2 U_j^M$, 
we have
\begin{align}
z^* - e^{-\mu_j^M}\left[ z^* - 0.5\left(1 - e^{-2U_j^L}\right)\right] 
\geq & z^* - e^{- 2U_j^M}\left[ z^* - 0.5\left(1 - e^{-2U_j^L}\right)\right]\nonumber \\
= & z^* - e^{-2(U_j^S - U_j^T)} \left[ z^* - 0.5\left(1 - e^{-2(U_j - U_j^S)}\right)\right]. \label{eq:newfinalproof1}
\end{align}

It is easy to see that (\ref{eq:newfinalproof1}) is decreasing in $U_j^T$. We next show that, given $U_j$ and $U_j^T$, (\ref{eq:newfinalproof1}) is also decreasing in $U_j^S$.
\begin{align*}
& \frac{\partial}{\partial U_j^S} \left[ z^* - e^{-2(U_j^S - U_j^T)} [ z^* - 0.5(1 - e^{-2(U_j - U_j^S)})] \right]\\
= & 2 e^{-2(U_j^S - U_j^T)} [ z^* - 0.5(1 - e^{-2(U_j - U_j^S)})] - e^{-2(U_j^S - U_j^T)} e^{-2(U_j - U_j^S)}\\
= & e^{-2(U_j^S - U_j^T)} ( 2 z^* - 1)\\
< & 0.
\end{align*}

According to Definition \ref{def:resourcetype}, we must have
\begin{equation}\label{eq:newboundproof1}
U_j^S < -0.5 \log(1 - 2r^* U_j) 
  \end{equation}
  and
  \begin{equation}\label{eq:newboundproof2}
U_j^T < - (1 - z^*) \log\left(1 - \frac{r^*U_j}{1-z^*}\right).
  \end{equation}

Since (\ref{eq:newfinalproof1}) is decreasing in $U_j^S$ and $U_j^T$, we can plug in (\ref{eq:newboundproof1}) and (\ref{eq:newboundproof2}) and obtain
\begin{align}
&z^* - e^{-2(U_j^S - U_j^T)} [ z^* - 0.5(1 - e^{-2(U_j - U_j^S)})]\nonumber \\
\geq & z^* - e^{-2(- 0.5  \log(1 - 2r^* U_j) + (1 - z^*) \log(1 - \frac{r^* U_j}{1-z^*}))} [z^* - 0.5(1 - e^{-2(U_j + 0.5  \log(1 - 2r^*U_j))})]\nonumber \\
= & z^* - \left[ z^* - \frac{1}{2} \left( 1 - \frac{1}{1-2r^*U_j}\cdot \frac{1}{e^{2U_j}}  \right) \right] (1 - 2r^* U_j) \left( \frac{1 - z^*}{1 - z^* - r^* U_j}\right) ^{2(1-z^*)}. \label{eq:newfinalproof1b}
\end{align}

Since $z^*$ and $r^*$ are constants, (\ref{eq:newfinalproof1b}) is a function of a single variable $U_j$. It is easy to check that this function is increasing and concave in $U_j$ for $U_j \leq c_j = 1$. Moreover, it equals $0$ at $U_j = 0$. Therefore, 
\begin{align*}
& z^* - \left[ z^* - \frac{1}{2} \left( 1 - \frac{1}{1-2r^*U_j}\cdot \frac{1}{e^{2U_j}}  \right) \right] (1 - 2r^* U_j) \left( \frac{1 - z^*}{1 - z^* - r^* U_j}\right) ^{2(1-z^*)}\\
\geq & U_j \left[  z^* - \left[ z^* - \frac{1}{2} \left( 1 -\frac{1}{1-2r^*}\cdot \frac{1}{e^2} \right) \right] (1 - 2r^*) \left( \frac{1 - z^*}{1 - z^* - r^*}\right) ^{2(1-z^*)} \right]\\
= & U_j h(z^*,r^*)\\
= & U_j r^*.
\end{align*}

This completes the proof for the theorem.
\halmos
\endproof

\section{Numerical Study}\label{sec:numerical}
We compare the empirical performance of our algorithms against two commonly used heuristics by simulating the algorithms on appointment-scheduling data obtained from a large hospital system in New York City. 

We obtain our data set from an Allergy department in the hospital system.  The data set contains more than 20000 appointment entries recorded in the year 2013. Each entry in the data records information about one appointment.  The entry includes the date that the patient makes the appointment, the exact time of the appointment, whether the patient eventually showed up to the original appointment, canceled the appointment some time later, or missed the appointment.

\begin{figure}[h!]
   \centering
   \caption{Average number of arrivals in a week.}
   \label{fig:arrivals}
     \includegraphics[width=0.5\textwidth]{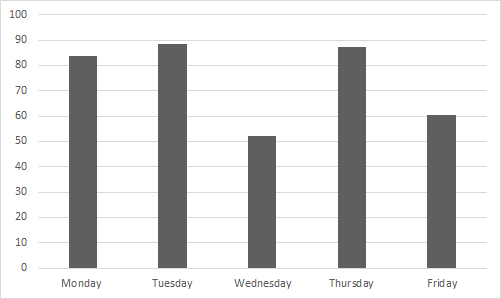}
\end{figure}

The average total number of patients who arrive to make appointments on each day is shown in Figure \ref{fig:arrivals}. It can be readily seen that the arrival pattern is highly non-stationary, as the average total number of arrivals on Thursday is $60\%$ more than that on Wednesday. 

We simulate a discrete horizon of 200 days. In each day, a random number of patients arrive to make appointments. Each patient needs to be assigned an appointment of 15 min, 30 min, 45 min, depending on his or her condition. By medical necessity, some patients must be assigned same-day appointments if at all.   We call these patients \emph{urgent patients}.  Other patients can be assigned to any day in the future.  We call these patients \emph{regular patients}. The relative proportions of patients in each priority category are summarized in Table \ref{tab:patientCategory}. We impose a requirement that regular patients must be assigned an appointment that is no more than 20-days away from the date of his or her first request for an appointment.  Although this hard deadline is not strictly enforced in reality, consideration for patient satisfaction often impels the administration to limit as much as possible the number of days that each patient must be made to wait.  Our deadline mimics this effect. 

\begin{table}
\begin{center}
\caption{Percentage of patients in different categories.}
\label{tab:patientCategory}
\begin{tabular}{|c|c|c|c|}
\hline
& 15 min & 30 min & 45 min\\
\hline
urgent  & $27\%$ & $1\%$ & $0\%$\\
\hline
regular & $45\%$ & $14\%$ & $ 9\%$\\
\hline
\end{tabular}
\end{center}
\end{table}

We assume a 5-day work week. We estimate the expected number of patients arriving per day of the week as shown in Figure \ref{fig:arrivals}. We assume that each patient randomly and independently falls into one of the six categories shown in Table \ref{tab:patientCategory}. \x{ All patients, whether urgent or regular, arrive at the beginning of a day. In our model,   the type of a patient is defined by both the time of arrival (one of these 200 days) and one of the six conventional ``types'' as defined in Table \ref{tab:patientCategory}. Overall, there are 1,200 patient types in our model. Moreover, each $\lambda_i(t)$ is only non-zero for one day of the horizon. }

We assume that there are multiple sessions on each day. Each session corresponds to a resource in our model.  We vary the session length among 1, 1.5, 2, 3, or 4 hours.  We assume that a patient can be assigned to any appointment within a day, as long as there is enough service time remaining and the day falls within the deadline to serve the patient.  We vary the number of sessions that are available per day.

We test the following two algorithms
\begin{itemize}
\item Our basic online algorithm ($LS$).
\item Our modified $LS$ algorithm ($MLS$).
\item Our refined  algorithm ($RLS$).
\item A greedy heuristic ($GRD$) that tries to assign every patient to the most recent session that is available and falls within his deadline.
\item A heuristic ($RSRV$) that reserves for each category an amount of capacity that is approximately equal to the average utilization of that category. This reservation is nested in the sense that higher-priority patients have access to their reserved capacity, as well as the reserved capacity of all lower-priority categories.  The heuristic then assigns patients greedily to the reserved capacity.  
\item The primal-dual algorithm ($PD$) given by \citet{buchbinder2007online}.
\end{itemize}

For each algorithm and each test case, we simulate the total length of appointments made during the entire 200 periods and calculate the average total length over 1000 replicates. We report the ratio of this average number relative to the optimal objective value of the upper bound given in \eqref{eq:lp}. \x{Note that in this numerical setting, the LP \eqref{eq:lp} can be solved by a simple greedy approach. First, we pack all urgent patient types into the same-day appointment sessions. Then, for period $t$ from $1$ to $200$, we pack regular patient types which arrive in period $t$ into the earliest available sessions. This yields an optimal solution to LP \eqref{eq:lp}. In more general settings for which simple heuristics do not give optimal LP solutions, one can always apply efficient packing LP solvers \citep{Allen-Zhu2018}. }

\begin{table}[h]
\begin{center}
\caption{Algorithm performance and average waiting time of regular patients. The length of each session is 1 hour. }
\label{tab:num1h}
\begin{tabular}{|c|c||c|c|c|c|c|}
\hline
Number of sessions & Scale & LS & RLS & GRD & RSRV & PD\\
\hline
18	&$	70.5\%	$&$	69.6\%	$, $	18.6	$&$	94.7\%	$, $	16.3	$&$	98.4\%	$, $	14.4	$&$	80.2\%	$, $	17.5	$&$	97.9\%	$, $	16.7	$\\
19	&$	74.4\%	$&$	69.2\%	$, $	18.3	$&$	94.2\%	$, $	15.4	$&$	98.4\%	$, $	11.0	$&$	78.8\%	$, $	17.4	$&$	97.6\%	$, $	16.9	$\\
20	&$	78.3\%	$&$	69.3\%	$, $	18.1	$&$	94.3\%	$, $	13.5	$&$	98.1\%	$, $	6.4	$&$	77.4\%	$, $	17.2	$&$	97.1\%	$, $	16.9	$\\
21	&$	82.2\%	$&$	69.4\%	$, $	17.7	$&$	94.4\%	$, $	11.1	$&$	96.2\%	$, $	2.4	$&$	76.2\%	$, $	16.9	$&$	96.4\%	$, $	16.7	$\\
22	&$	86.2\%	$&$	69.9\%	$, $	17.2	$&$	94.5\%	$, $	10.8	$&$	96.0\%	$, $	1.3	$&$	75.1\%	$, $	16.6	$&$	95.5\%	$, $	16.2	$\\
23	&$	90.1\%	$&$	70.3\%	$, $	16.2	$&$	94.5\%	$, $	11.0	$&$	96.0\%	$, $	0.9	$&$	74.0\%	$, $	16.2	$&$	94.3\%	$, $	15.5	$\\
24	&$	94.0\%	$&$	70.8\%	$, $	13.8	$&$	94.2\%	$, $	10.2	$&$	95.5\%	$, $	0.7	$&$	73.1\%	$, $	15.6	$&$	92.5\%	$, $	15.6	$\\
25	&$	97.9\%	$&$	70.7\%	$, $	9.0	$&$	93.7\%	$, $	6.8	$&$	94.5\%	$, $	0.5	$&$	72.2\%	$, $	14.6	$&$	90.5\%	$, $	15.9	$\\
26	&$	101.8\%	$&$	70.7\%	$, $	4.2	$&$	93.4\%	$, $	3.3	$&$	94.7\%	$, $	0.4	$&$	73.0\%	$, $	12.8	$&$	90.3\%	$, $	16.1	$\\
27	&$	105.7\%	$&$	71.0\%	$, $	1.6	$&$	95.3\%	$, $	1.3	$&$	95.7\%	$, $	0.3	$&$	74.8\%	$, $	9.7	$&$	91.3\%	$, $	16.0	$\\
28	&$	109.7\%	$&$	71.0\%	$, $	0.9	$&$	96.6\%	$, $	0.8	$&$	96.6\%	$, $	0.2	$&$	76.6\%	$, $	6.3	$&$	92.4\%	$, $	15.7	$\\
29	&$	113.6\%	$&$	70.9\%	$, $	0.6	$&$	97.4\%	$, $	0.6	$&$	97.2\%	$, $	0.2	$&$	77.4\%	$, $	3.2	$&$	93.4\%	$, $	15.4	$\\
30	&$	117.5\%	$&$	71.0\%	$, $	0.5	$&$	97.8\%	$, $	0.4	$&$	97.7\%	$, $	0.2	$&$	77.8\%	$, $	1.4	$&$	94.2\%	$, $	15.1	$\\
31	&$	121.4\%	$&$	71.0\%	$, $	0.4	$&$	98.3\%	$, $	0.3	$&$	98.2\%	$, $	0.1	$&$	78.0\%	$, $	0.9	$&$	95.1\%	$, $	14.4	$\\
32	&$	125.3\%	$&$	71.0\%	$, $	0.3	$&$	98.6\%	$, $	0.3	$&$	98.4\%	$, $	0.1	$&$	78.1\%	$, $	0.6	$&$	96.4\%	$, $	12.9	$\\
33	&$	129.2\%	$&$	70.9\%	$, $	0.2	$&$	98.8\%	$, $	0.2	$&$	98.7\%	$, $	0.1	$&$	78.1\%	$, $	0.4	$&$	97.2\%	$, $	12.0	$\\
\hline
\end{tabular}
\end{center}
\end{table}

\begin{table}[h]
\begin{center}
\caption{Algorithm performance and average waiting time of regular patients. The length of each session is 1.5 hours. }
\label{tab:num15h}
\begin{tabular}{|c|c||c|c|c|c|c|}
\hline
Number of sessions & Scale & MLS & RLS & GRD & RSRV & PD\\
\hline
12	&$	70.5\%	$&$	76.2\%	$, $	18.6	$&$	98.2\%	$, $	15.9	$&$	98.4\%	$, $	14.4	$&$	92.7\%	$, $	17.7	$&$	97.9\%	$, $	17.5	$\\
13	&$	76.4\%	$&$	76.3\%	$, $	18.3	$&$	97.9\%	$, $	12.8	$&$	98.4\%	$, $	8.8	$&$	91.5\%	$, $	17.5	$&$	97.4\%	$, $	17.5	$\\
14	&$	82.2\%	$&$	76.5\%	$, $	17.7	$&$	97.7\%	$, $	10.8	$&$	96.3\%	$, $	2.4	$&$	90.2\%	$, $	17.3	$&$	96.7\%	$, $	17.4	$\\
15	&$	88.1\%	$&$	77.1\%	$, $	16.8	$&$	97.5\%	$, $	11.5	$&$	96.1\%	$, $	1.1	$&$	88.9\%	$, $	17.0	$&$	95.7\%	$, $	17.3	$\\
16	&$	94.0\%	$&$	77.5\%	$, $	13.9	$&$	96.7\%	$, $	10.8	$&$	95.6\%	$, $	0.7	$&$	87.6\%	$, $	16.5	$&$	93.9\%	$, $	17.2	$\\
17	&$	99.9\%	$&$	77.5\%	$, $	6.8	$&$	95.7\%	$, $	5.5	$&$	93.7\%	$, $	0.5	$&$	86.5\%	$, $	15.7	$&$	91.5\%	$, $	17.1	$\\
18	&$	105.7\%	$&$	77.6\%	$, $	1.7	$&$	97.2\%	$, $	1.4	$&$	95.7\%	$, $	0.3	$&$	90.3\%	$, $	14.3	$&$	93.7\%	$, $	16.8	$\\
19	&$	111.6\%	$&$	77.7\%	$, $	0.8	$&$	98.0\%	$, $	0.7	$&$	96.8\%	$, $	0.2	$&$	94.2\%	$, $	11.0	$&$	95.0\%	$, $	16.7	$\\
20	&$	117.5\%	$&$	77.7\%	$, $	0.5	$&$	98.4\%	$, $	0.5	$&$	97.6\%	$, $	0.2	$&$	97.9\%	$, $	6.3	$&$	95.9\%	$, $	16.5	$\\
21	&$	123.4\%	$&$	77.7\%	$, $	0.4	$&$	98.7\%	$, $	0.3	$&$	98.1\%	$, $	0.1	$&$	99.2\%	$, $	2.2	$&$	96.5\%	$, $	16.5	$\\
22	&$	129.2\%	$&$	77.7\%	$, $	0.3	$&$	99.0\%	$, $	0.3	$&$	98.5\%	$, $	0.1	$&$	99.6\%	$, $	0.9	$&$	97.1\%	$, $	16.4	$\\
\hline
\end{tabular}
\end{center}
\end{table}

\begin{table}[h]
\begin{center}
\caption{Algorithm performance and average waiting time of regular patients. The length of each session is 2 hours. }
\label{tab:num2h}
\begin{tabular}{|c|c||c|c|c|c|c|}
\hline
Number of sessions & Scale & MLS & RLS & GRD & RSRV & PD\\
\hline
9	&$	70.5\%	$&$	78.0\%	$, $	18.6	$&$	99.1\%	$, $	16.0	$&$	98.5\%	$, $	14.4	$&$	91.1\%	$, $	17.9	$&$	97.9\%	$, $	17.7	$\\
10	&$	78.3\%	$&$	77.3\%	$, $	18.1	$&$	99.0\%	$, $	11.2	$&$	98.4\%	$, $	6.4	$&$	89.2\%	$, $	17.7	$&$	97.3\%	$, $	17.6	$\\
11	&$	86.2\%	$&$	77.9\%	$, $	17.2	$&$	98.6\%	$, $	11.6	$&$	96.2\%	$, $	1.4	$&$	87.2\%	$, $	17.4	$&$	96.3\%	$, $	17.5	$\\
12	&$	94.0\%	$&$	78.3\%	$, $	13.9	$&$	97.3\%	$, $	11.0	$&$	95.7\%	$, $	0.7	$&$	85.4\%	$, $	17.1	$&$	94.5\%	$, $	17.4	$\\
13	&$	101.8\%	$&$	78.3\%	$, $	4.4	$&$	96.4\%	$, $	3.7	$&$	94.4\%	$, $	0.4	$&$	85.3\%	$, $	16.4	$&$	93.0\%	$, $	17.3	$\\
14	&$	109.7\%	$&$	78.9\%	$, $	1.0	$&$	98.0\%	$, $	0.9	$&$	96.5\%	$, $	0.3	$&$	90.2\%	$, $	15.0	$&$	95.1\%	$, $	17.2	$\\
15	&$	117.5\%	$&$	79.2\%	$, $	0.5	$&$	98.6\%	$, $	0.5	$&$	97.6\%	$, $	0.2	$&$	95.2\%	$, $	10.8	$&$	96.3\%	$, $	17.1	$\\
16	&$	125.3\%	$&$	79.1\%	$, $	0.3	$&$	99.0\%	$, $	0.3	$&$	98.2\%	$, $	0.1	$&$	99.4\%	$, $	4.8	$&$	97.0\%	$, $	16.9	$\\
\hline
\end{tabular}
\end{center}
\end{table}

\begin{table}[h]
\begin{center}
\caption{Algorithm performance and average waiting time of regular patients.  The length of each session is 3 hours.  }
\label{tab:num3h}
\begin{tabular}{|c|c||c|c|c|c|c|}
\hline
Number of sessions & Scale & MLS & RLS & GRD & RSRV & PD\\
\hline
6	&$	70.5\%	$&$	84.3\%	$, $	18.6	$&$	99.3\%	$, $	16.3	$&$	98.6\%	$, $	14.3	$&$	93.9\%	$, $	17.7	$&$	98.0\%	$, $	18.0	$\\
7	&$	82.2\%	$&$	84.5\%	$, $	17.8	$&$	98.9\%	$, $	11.6	$&$	96.5\%	$, $	2.4	$&$	91.8\%	$, $	17.2	$&$	96.9\%	$, $	17.7	$\\
8	&$	94.0\%	$&$	84.9\%	$, $	13.9	$&$	97.4\%	$, $	11.6	$&$	95.7\%	$, $	0.7	$&$	89.5\%	$, $	16.3	$&$	94.8\%	$, $	17.6	$\\
9	&$	105.7\%	$&$	85.3\%	$, $	1.7	$&$	97.5\%	$, $	1.5	$&$	95.6\%	$, $	0.3	$&$	92.4\%	$, $	13.3	$&$	94.7\%	$, $	17.5	$\\
10	&$	117.5\%	$&$	85.4\%	$, $	0.5	$&$	98.8\%	$, $	0.5	$&$	97.6\%	$, $	0.2	$&$	99.6\%	$, $	4.3	$&$	96.5\%	$, $	17.4	$\\
11	&$	129.2\%	$&$	85.4\%	$, $	0.3	$&$	99.4\%	$, $	0.3	$&$	98.5\%	$, $	0.1	$&$	99.9\%	$, $	0.7	$&$	97.4\%	$, $	17.2	$\\
\hline
\end{tabular}
\end{center}
\end{table}

\begin{table}[h]
\begin{center}
\caption{Algorithm performance and average waiting time of regular patients.  The length of each session is 4 hours.  }
\label{tab:num4h}
\begin{tabular}{|c|c||c|c|c|c|c|}
\hline
Number of sessions & Scale & MLS & RLS & GRD & RSRV & PD\\
\hline
5	&$	78.3\%	$&$	86.0\%	$, $	18.1	$&$	99.2\%	$, $	11.9	$&$	98.6\%	$, $	6.2	$&$	93.7\%	$, $	17.4	$&$	97.4\%	$, $	17.8	$\\
6	&$	94.0\%	$&$	86.7\%	$, $	13.9	$&$	97.3\%	$, $	11.9	$&$	95.8\%	$, $	0.7	$&$	91.1\%	$, $	16.0	$&$	94.9\%	$, $	17.6	$\\
7	&$	109.7\%	$&$	87.1\%	$, $	1.0	$&$	98.2\%	$, $	0.9	$&$	96.5\%	$, $	0.3	$&$	97.0\%	$, $	8.2	$&$	95.6\%	$, $	17.5	$\\
8	&$	125.3\%	$&$	87.5\%	$, $	0.3	$&$	99.3\%	$, $	0.3	$&$	98.3\%	$, $	0.1	$&$	99.8\%	$, $	0.8	$&$	97.1\%	$, $	17.3	$\\
\hline
\end{tabular}
\end{center}
\end{table}

\begin{table}[h]
\begin{center}
\caption{Algorithm performance and average waiting time of regular patients.  Regular patients arrive only on Mondays, and same-day patients arrive only on the other weekdays.  The length of each session is 1 hour.  }
\label{tab:numAlternating}
\begin{tabular}{|c|c||c|c|c|c|c|}
\hline
Number of sessions & Scale & LS & RLS & GRD & RSRV & PD\\
\hline
16	&$	74.4\%	$&$	75.5\%	$, $	17.4	$&$	94.7\%	$, $	6.3	$&$	93.9\%	$, $	0.6	$&$	86.1\%	$, $	1.1	$&$	86.7\%	$, $	7.3	$\\
17	&$	79.1\%	$&$	75.7\%	$, $	17.5	$&$	94.3\%	$, $	8.1	$&$	92.2\%	$, $	0.6	$&$	85.3\%	$, $	1.0	$&$	85.3\%	$, $	3.2	$\\
18	&$	83.7\%	$&$	76.0\%	$, $	17.2	$&$	94.2\%	$, $	9.8	$&$	90.6\%	$, $	0.5	$&$	84.0\%	$, $	0.9	$&$	84.6\%	$, $	2.1	$\\
19	&$	88.4\%	$&$	76.1\%	$, $	16.3	$&$	94.0\%	$, $	11.0	$&$	89.0\%	$, $	0.5	$&$	82.3\%	$, $	0.9	$&$	83.7\%	$, $	1.9	$\\
20	&$	93.0\%	$&$	76.1\%	$, $	11.1	$&$	93.9\%	$, $	8.4	$&$	88.0\%	$, $	0.5	$&$	80.8\%	$, $	0.8	$&$	83.0\%	$, $	1.7	$\\
21	&$	97.7\%	$&$	76.2\%	$, $	3.5	$&$	94.4\%	$, $	2.7	$&$	88.6\%	$, $	0.4	$&$	80.6\%	$, $	0.8	$&$	83.6\%	$, $	1.6	$\\
22	&$	102.3\%	$&$	76.1\%	$, $	2.0	$&$	94.7\%	$, $	1.7	$&$	89.5\%	$, $	0.4	$&$	80.7\%	$, $	0.7	$&$	84.4\%	$, $	1.5	$\\
23	&$	107.0\%	$&$	75.7\%	$, $	1.5	$&$	94.5\%	$, $	1.3	$&$	90.1\%	$, $	0.4	$&$	80.8\%	$, $	0.7	$&$	85.2\%	$, $	1.4	$\\
24	&$	111.6\%	$&$	75.6\%	$, $	1.2	$&$	95.0\%	$, $	1.0	$&$	91.9\%	$, $	0.3	$&$	82.0\%	$, $	0.6	$&$	87.0\%	$, $	1.3	$\\
25	&$	116.3\%	$&$	75.6\%	$, $	1.0	$&$	95.7\%	$, $	0.8	$&$	93.8\%	$, $	0.3	$&$	83.4\%	$, $	0.6	$&$	88.9\%	$, $	1.3	$\\
26	&$	120.9\%	$&$	75.6\%	$, $	0.8	$&$	96.2\%	$, $	0.7	$&$	95.4\%	$, $	0.3	$&$	84.8\%	$, $	0.6	$&$	90.5\%	$, $	1.2	$\\
27	&$	125.6\%	$&$	75.6\%	$, $	0.6	$&$	96.6\%	$, $	0.6	$&$	96.7\%	$, $	0.2	$&$	86.2\%	$, $	0.5	$&$	91.8\%	$, $	1.2	$\\
\hline
\end{tabular}
\end{center}
\end{table}

\begin{table}[h]
\begin{center}
\caption{Performance relative to the upper bound given in \eqref{eq:lp}, when parameters are randomly generated.}
\label{tab:numrand}
\begin{tabular}{|c||c|c|c|c|c|}
\hline
Number of sessions &  LS & RLS & GRD & RSRV & PD\\
\hline
Worst Setting	&$		44.3\%	$&$	68.4\%	$&$	66.3\%	$&$	43.7\%	$&$	67.6\%	$\\
Average Setting&$		65.2\%	$&$	96.3\%	$&$	95.9\%	$&$	85.2\%	$&$	95.9\%	$\\
\hline
\end{tabular}
\end{center}
\end{table}

Tables \ref{tab:num1h} to \ref{tab:num4h} summarize the performance of the algorithms.  The \emph{scale} is the ratio of total capacity to total demand. \x{In each cell, the first number is the performance of the algorithm relative to the upper bound \eqref{eq:lp}; the second number is the average number of days that admitted regular patients need to wait under the algorithm.}  We make several observations:
\begin{itemize}
\item The refined algorithm $RLS$ is never more than 7\% worse than the upper bound on average in each of the scenarios tested.   The reservation heuristic $RSRV$ could be as much as 16\% worse than the upper bound on average.  The greedy heuristic $GRD$ could be as much as 5.7\% worse than the upper bound on average.

\item Predictably, the refined algorithm $RLS$ dominates the basic algorithms $LS$ and $MLS$.  This performance gain comes from better resource sharing. 

\item The greedy heuristic $GRD$ also performs consistently better than the static reservation heuristic $RSRV$, except when the scale is high.  Most likely, the greedy heuristic allows greater resource sharing among different customer types, which results in better resource utilization.  However, when the scale is high, there is an abundance of capacity, so that resource sharing is less important.  

\item The greedy heuristic $GRD$ tends to be good when the scale is either very large or very small.   These are situations in which it is easier to do well.  When there is little capacity, the utilization can be kept high even with a naive algorithm because there is relative very high demand.  When there is an abundance of capacity, the utilization can be close to optimal because a high proportion of demand can be accommodated.   Therefore, an algorithm offers the most value relative to a naive heuristic when the scale is moderate.  

\item Similar to $GRD$, the Primal-Dual algorithm performs well when the scale is either large or small.  However, its performance is slightly worse than $GRD$ in most cases. This might be because the Primal-Dual algorithm is specially designed to improve the worst-case performance, whereas we report the average-case performance.

\item The refined algorithm $RLS$ performs significantly better than, or is very close to, the better of the two heuristics.  It performs much better than the heuristics when the scale is moderate, which is when an algorithm offers the most value relative to a naive heuristic.

\item \x{The average number of days that admitted regular patients need to wait under the greedy heuristic is the smallest among all algorithms.  This is because the greedy heuristic allows regular patients to take sessions that could have been reserved for same-day patients. In Table 8, we test scenarios that better illustrate the outcomes of making greedy assignments. We let regular patients arrive only on Mondays, and let same-day patients arrive only on the other four weekdays. We find that the greedy heuristic results in extremely short waiting times compared to other algorithms, but the performance of the greedy heuristic is much worse due to the fact that it does not reserve the right amount of resources for same-day patients.} 
\end{itemize}

We also test the algorithms under randomly generated settings. In Table \ref{tab:numrand}, we report the worst performance and the average performance of all the algorithms over 100 random settings. The performance of algorithms in each setting is calculated by simulating 1000 replicates. Each of the 100 random settings is generated by 
\begin{itemize}
\item uniformly generating the percentages in Table \ref{tab:patientCategory};
\item uniformly picking a deadline for all regular patients between 5 and 30 days;
\item uniformly setting the capacity of all resources to be between 45 and 150 minutes;
\item uniformly picking a scale between $70\%$ and $130\%$.
\end{itemize}

Again, our RLS algorithm consistently performs well in these test cases.


%
%
%

\ACKNOWLEDGMENT{The second two authors gratefully acknowledge support by the National Science Foundation under award CMMI 1538088.}





\bibliographystyle{ormsv080}
\bibliography{myrefs}{}
\end{document}